\documentclass[12pt,leqno]{compositio}
\usepackage{amsmath,amssymb,amscd,latexsym}

\usepackage[all]{xy}
\newcommand{\cart}[3]{\xymatrix{ #1 \ar@{}[r]|{\times}_-{#2} & #3}}
\newcommand{\cartt}[5]{\xymatrix{ #1 \ar@{}[r]|{\times}_-{#2} & #3 \ar@{}[r]|{\times}_-{#4} & #5}}
\newcommand{\ocart}[3]{\xymatrix{ #1 \ar@{}[r]|{\otimes}_-{#2} & #3}}

\newcommand{\CH}{{\check{H}}}
\newcommand{\Aone}{{\hat{\mathbb{A}}^1}}

\newcommand{\BP}{{\mathrm{BP}}}
\newcommand{\FGL}{{\mathrm{FGL}}}
\newcommand{\SI}{{\mathrm{SI}}}
\newcommand{\K}{{\mathrm{K}}}
\newcommand{\MU}{{\mathrm{MU}}}
\newcommand{\Rings}{{\mathrm{Rings}}}
\newcommand{\FRings}{{\mathrm{FRings}}}
\newcommand{\op}{{\mathrm{op}}}
\newcommand{\Spf}{{\mathrm{Spf}}}
\newcommand{\Pic}{{\mathrm{Pic}}}
\newcommand{\F}{{\mathbb{F}}}
\newcommand{\Ge}{\mathbb{G}}
\newcommand{\HH}{{\mathrm{H}}}

\newcommand{\Z}{{\mathbb{Z}}}
\newcommand{\Aff}{\mathrm{Aff}}
\newcommand{\Sets}{\mathrm{Sets}}
\newcommand{\Ell}{\mathrm{Ell}}
\newcommand{\Ello}{\overline{\mathrm{Ell}}}

\newcommand{\Isom}{\underline{\mathrm{Isom}}}
\newcommand{\height}{\mathrm{ht}}
\newcommand{\id}{\mathrm{id}}
\renewcommand{\ker}{\mathrm{ker}}
\renewcommand{\max}{\mathrm{max}}
\renewcommand{\min}{\mathrm{min}}

\newcommand{\Modq}[1]{\mathrm{Mod}_{\mathrm{qcoh}}({\cal O}_{#1})}
\newcommand{\ModX}{\mathrm{Mod}_{\mathrm{qcoh}}({\cal O}_{\mathfrak{X}})}

\newcommand{\ModXz}{\mathrm{Mod}_{\mathrm{qcoh}}({\cal O}_{X_0})}

\newcommand{\ModXo}{\mathrm{Mod}_{\mathrm{qcoh}}({\cal O}_{X_1})}

\newcommand{\Mor}{\mathrm{Mor}\;}
\newcommand{\Ob}{\mathrm{Ob}}
\newcommand{\Spec}{\mathrm{Spec}\,}
\newcommand{\Aut}{\mathrm{Aut}\,}
\newcommand{\Ext}{\mathrm{Ext}}
\newcommand{\Hom}{\mathrm{Hom}}
\newcommand{\Fh}{{\mathcal F}}
\newcommand{\Gh}{{\mathcal G}}

\newcommand{\Ih}{{\mathcal I}}
\newcommand{\Oh}{{\mathcal O}}

\newcommand{\Ug}{{\mathfrak{U}}}
\newcommand{\Yg}{{\mathfrak{Y}}}
\newcommand{\Xg}{{\mathfrak{X}}}
\newcommand{\Zg}{{\mathfrak{Z}}}

\newtheorem{theorem}{Theorem}

\newtheorem{prop}[theorem]{Proposition}
\newtheorem{defn}[theorem]{Definition}
\newtheorem{cor}[theorem]{Corollary}

\newtheorem{remark}[theorem]{Remark}

\newenvironment{proofof}{\noindent {\bf Proof of}}{\mbox{}\hfill$\Box$}
\begin{document}

\title[The stack of formal groups]{The stack of formal groups in stable homotopy theory}
\author{N. Naumann}
\email{niko.naumann@mathematik.uni-regensburg.de}
\address{NWF I- Mathematik\\ Universit\"at Regensburg\\ Universit\"atsstrasse 31\\93053 Regensburg}

\begin{abstract}
We construct the algebraic stack of formal groups and use it to provide
a new perspective onto a recent result of M. Hovey and N. Strickland on comodule categories for Landweber exact algebras. This leads to a geometric understanding of their results as well as to a generalisation.
\end{abstract}

\maketitle

\section{Introduction}

Ever since the fundamental work of S. Novikov and D. Quillen \cite{novikov},\cite{Q} the theory of formal groups
is firmly rooted in stable homotopy theory. In particular, the simple geometric
structure of the moduli space of formal groups has been a constant source of inspiration. This moduli space is stratified according to the height of
the formal group. For many spaces $X$, $\MU_*(X)$ can canonically be considered
as a flat sheaf on the moduli space and the stratification defines a resolution
of $\MU_*(X)$ (the Cousin-complex) which is well-known to be the chromatic resolution
of $\MU_*(X)$ and which is a central tool in the actual computation of the stable
homotopy of $X$.\\
In fact, much deeper homotopy theoretic results have been 
suggested by this point of view and we mention two of them. All thick 
subcategories of the derived category of sheaves on the moduli space
are rather easily determined by using the above stratification. This
simple structure persists to determine all the thick subcategories of
the category of finite spectra, see \cite{R2}, Theorem 3.4.3. Similarly,
any coherent sheaf on the moduli space can be reconstructed from
its restriction to the various strata (corresponding roughly to $\K(n)$-
localisation in stable homotopy). Again, this result persists
to homotopy theory as the chromatic convergence theorem, \cite{R2}, Theorem 7.5.7.\\
In conclusion, the derived category of sheaves on the moduli space of 
formal groups has turned out to be an excellent algebraic approximation to the 
homotopy category of (finite) spectra.\\
It may thus seem a little surprising that the central notion of a stack of 
formal groups has not yet been given a solid foundation, and the chief
purpose of this paper is to do so. We hasten to point out to the knowledgeable
reader that to this end there is something to do beyond just copying 
existing literature as already the following simple remark demonstrates.
Defining (as usual) a formal group to be a group structure on the formal 
affine line, one is guaranteed to {\em not} obtain a stack just because the formal
affine line in general does admit non-trivial flat forms. We thus spend some
effort in the construction of the stack of formal groups and the derivation
of its basic properties. This may also be useful in the multiplicative
ring spectrum project of P. Goerss and M. Hopkins, c.f. \cite{G}.\\
In fact, we start out more generally by making precise the relation 
between flat Hopf algebroids and a certain class of stacks. Roughly, the datum of a flat Hopf algebroid is equivalent to the datum of the stack with a specific presentation. Now, the category of
comodules of the flat Hopf algebroid only depends on the stack. We will
demonstrate the gain in conceptual clarity provided by this point of
view by reconsidering the following
remarkable recent result of M. Hovey and N. Strickland. For two Landweber
exact $\BP_*$-algebras $R$ and $S$
of the same height the categories of comodules of the flat Hopf algebroids
$(R,\Gamma_R:=R\otimes_{\BP_*}\BP_*\BP\otimes_{\BP_*} R)$ and $(S,\Gamma_S:=S\otimes_{\BP_*}\BP_*\BP\otimes_{\BP_*} S)$ are equivalent. As an immediate consequence one obtains the computationally important
change-of-rings isomorphism $\Ext_{\Gamma_R}^*(R,R)\simeq\Ext_{\Gamma_S}^*(S,S)$
which had been established previously by G. Laures \cite{coco}, 4.3.3.\\
From our point of view, this result has the following simple explanation.
Let $\Xg$ be the stack associated with $(\BP_*,\BP_*\BP)$ and $f:\Spec(R)\longrightarrow\Xg$
the canonical map. As we will explain, $\Xg$ is closely related to the stack of formal
groups and is thus stratified by closed substacks
\[ \Xg=\Zg^0\supseteq\Zg^1\supseteq\ldots \,.\]
We will show that the induced Hopf algebroid $(R,\Gamma_R)$ is simply a presentation
of the {\em stack}-theoretic image of $f$ and that $R$ being Landweber
exact of height $n$ implies that this image is $\Xg-\Zg^{n+1}$. We conclude that
$(R,\Gamma_R)$ and $(S,\Gamma_S)$ are presentations of {\em the same}
stack which implies the result of \cite{H1} but more is true: The 
comodule categories under consideration are in fact equivalent as {\em tensor}
abelian categories (\cite{H1} treats their structure of abelian categories only)
and we easily generalise the above proof to apply to all the
stacks $\Zg^{n}-\Zg^{n+k}$ (with $n\geq 1$ allowed).\\
Returning to the stack of formal groups, we show that the stack associated with
$(\MU_*,\MU_*\MU)$ is closely related to this stack. Note, however, that 
this requires an {\em a priori} construction of the stack of formal groups, the
problem being the following. The objects of a stack associated with a 
flat Hopf algebroid are only {\em flat locally} given in terms of the 
Hopf algebroid and it is in general difficult to decide what additional 
objects the stack contains. Given the central role of the stack of formal groups in 
stable homotopy theory, we believe that it is important to have a genuinely
geometric understanding of it rather than just as the stack associated 
to some Hopf algebroid, so we solve this problem here.\\
Finally, we point out that for many stacks appearing in algebraic topology
it is surprisingly easy to compute their Picard groups. For example,
the Picard group of the stack of formal groups is
isomorphic to $\Z$, generated by the canonical line bundle. In keeping
with the philosophy that the stack of formal groups provides a good algebraic
approximation to stable homotopy theory one may try to use this in the current
investigations of the Picard groups of various categories of spectra (e.g. \cite{invertible}) and we hope to return to this in the future.\\
We review the individual sections in more detail. In section \ref{prelim}
we review the stack theoretic notions we will have to use in the following. In section \ref{stacksandhopf}
we give the relation between flat Hopf algebroids and algebraic stacks. In section \ref{morphisms} we collect a number
of technical results on algebraic stacks centring around the 
problem to relate the properties of a morphism between
algebraic stacks with properties of the functors
it induces on the categories of quasi-coherent
sheaves. The main result is proved in section \ref{ringchange}. In the final section \ref{stackformalgroups} we construct 
the stack of formal groups and show that the algebraic stack associated with the flat Hopf
algebroid $(\MU_*,\MU_*\MU)$ is the stack of (one dimensional, commutative, connected, formally smooth) formal groups together with a trivialization of the 
canonical line bundle and explain its basic geometric properties.\\
To conclude the introduction we would like to acknowledge the influence of M. Hopkins
on the present circle of ideas. We understand that he was the first to insist that numerous results 
on (comodules over) flat Hopf algebroids should be understood from a geometric, i.e. stack
theoretic, point of view, c.f. \cite{hopkinslecture}.\\

\begin{acknowledgements} I would like
to thank E. Pribble for making a copy of \cite{P} available to me,
J. Heinloth, J. Hornbostel, M. Hovey, G. Kings and N. Strickland for useful conversation and the referee for 
suggesting substantial improvements of the exposition.
\end{acknowledgements}

\section{Preliminaries on algebraic stacks}\label{prelim}

In this section we will recall those concepts from the theory of stacks which will be
used in the sequel.\\
Fix an affine scheme $S$ and denote by $\Aff_S$ the category of affine $S$-schemes with some
cardinality bound to make it small. We may write $\Aff$ for $\Aff_S$ if $S$ is understood.

\begin{defn}\label{def1}
A category fibred in groupoids (understood: over $\Aff$) is a category $\Xg$ together with a functor
$a:\Xg\longrightarrow\Aff$ such that\\
i) ("existence of pull-backs") For every morphism $\phi:V\longrightarrow U$ in $\Aff$ and 
$x\in\Ob(\Xg)$ with $a(x)=U$ there is a morphism $f:y\longrightarrow x$ with $a(f)=\phi$.\\
ii) ("uniqueness of pull-backs up to unique isomorphism") For every diagram in $\Xg$
\[ \xymatrix{& z \ar[d]^h \\ y \ar[r]^f & x} \] lying via $a$ over a diagram
\[ \xymatrix{& W \ar[d]^{\chi} \ar[dl]_{\psi}\\ V \ar[r]^{\phi} & U} \] in $\Aff$ there is a unique
morphism $g:z\longrightarrow y$ in $\Xg$ such that $f\circ g=h$ and $a(g)=\psi$.
\end{defn}

As an example, consider the category $\Ell$ of elliptic curves having objects $E/U$ consisting of
an affine $S$-scheme $U$ and an elliptic curve $E$ over $U$. Morphisms in $\Ell$ are
cartesian diagrams
\begin{equation}\label{z1}
\xymatrix{E' \ar[r] \ar[d] & E \ar[d] \\ U' \ar[r]^f & U,}
\end{equation}

equivalently isomorphisms of elliptic curves over $U'$ from $E'$ to $E\times_U U'$. For
an explicit account of $\Aut_{\Ell}(E/U)$ see \cite{strickell}, section 5.\\
There is a functor \[ a:\Ell\longrightarrow\Aff\]
sending $E/U$ to $U$ and a morphism in $\Ell$ as in (\ref{z1}) to $f$.\\
Checking that $a$ makes $\Ell$ a category fibred in groupoids reveals that the main subtlety 
in Definition \ref{def1} lies in then non-uniqueness of cartesian products. A similar
example can be given using vector bundles on topological spaces \cite{hollander}, Example B.2.\\
Let $a:\Xg\longrightarrow\Aff$ be a category fibred in groupoids. For $U\in\Ob(\Aff)$ the 
fibre category $\Xg_U\subseteq\Xg$ is defined as the subcategory having objects $x\in\Ob(\Xg)$
with $a(x)=U$ and morphisms $f\in\Mor(\Xg)$ with $a(f)=\id_U$. The category $\Xg_U$ is a groupoid.
Choosing a pull-back as in Definition \ref{def1}, i) for every $\phi:V\longrightarrow U$
in $\Aff$ one can define functors $\phi^*:\Xg_U\longrightarrow\Xg_V$ and, for composable
$\phi,\psi\in\Mor(\Aff)$, isomorphisms $\psi^*\circ\phi^*\simeq (\phi\circ\psi)^*$ satisfying a cocycle condition. Sometimes $\phi^*(x)$ will be denoted as $x|V$. This connects Definition
\ref{def1} with the concept of fibred category as in \cite{SGA1}, VI as well as with the notion
of lax/pseudo functor/presheaf on $\Aff$ with values in groupoids; see \cite{hollander}
and \cite{vistoli} for more details.\\
Categories fibred in groupoids constitute a $2$-category in which $1$-morphisms from 
$a:\Xg\longrightarrow\Aff$ to $b:\Yg\longrightarrow\Aff$ are functors $f:\Xg\longrightarrow\Yg$
with $b\circ f=a$ (sic !) and $2$-morphisms are isomorphisms between $1$-morphisms. A $1$-morphism $f:\Xg\longrightarrow\Yg$
is called a monomorphism (resp. isomorphism) if for all $U\in\Ob(\Aff)$ the induced functor
$f_U:\Xg_U\longrightarrow\Yg_U$ between fibre categories is fully faithful (resp. an equivalence
of categories).\\
The next point is to explain what a sheaf, rather than a presheaf, of groupoids should be.
This makes sense for any topology on $\Aff$ but we fix the $fpqc$ topology for definiteness:
It is the Grothendieck topology on $\Aff$ generated by the pretopology which as covers of an 
$U\in\Aff$ has the finite families of flat morphisms $U_i\longrightarrow U$ in $\Aff$ such
that $\coprod_i U_i\longrightarrow U$ is faithfully flat, c.f. \cite{vistoli}, 2.3.

\begin{defn}\label{def2}
A stack (understood: over $\Aff$ for the $fpqc$ topology) is a category fibred in groupoids $\Xg$
such that\\
i) ("descent of morphisms") For $U\in\Ob(\Aff)$ and $x,y\in\Ob(\Xg_U)$ the presheaf
\[ \Aff/U\longrightarrow\Sets\, , \, (V\stackrel{\phi}{\longrightarrow}U)\mapsto \Hom_{\Xg_V}(x|V,y|V) \]
is a sheaf.\\
ii) ("glueing of objects") If $\{ U_i \stackrel{\phi_i}{\longrightarrow}U\}$ is a covering in
$\Aff$, $x_i\in\Ob(\Xg_{U_i})$ and $f_{ji}:(x_i|U_i\times_U U_j)\stackrel{\sim}{\longrightarrow}
(x_j|U_i\times_U U_j)$ are isomorphisms satisfying a cocycle condition then there are $x\in\Ob(\Xg_U)$ and isomorphisms $f_i:(x|U_i)\stackrel{\sim}{\longrightarrow}x_i$ such that $f_j|U_i\times_U U_j=f_{ji}\circ f_i|U_i\times_U U_j$.
\end{defn}

The category fibred in groupoids $\Ell$ is a stack: Condition i) of Definition \ref{def2}
for $\Ell$ is a consequence of faithfully flat descent \cite{BLR}, 6.1, Theorem 6, and
condition ii) relies on the fact that elliptic curves canonically admit ample line bundles, see
\cite{vistoli}, 4.3.3.\\

\begin{defn}\label{def3}
Let $\Xg$ be a stack. A substack of $\Xg$ is a strictly full subcategory $\Yg\subseteq\Xg$ such that\\
i) For any $\phi:U\longrightarrow V$ in $\Aff$ one has $\phi^*(\Ob(\Yg_V))\subseteq\Ob(\Yg_U)$.\\
ii) If $\{ U_i\longrightarrow U\}$ is a covering in $\Aff$ and $x\in\Ob(\Xg_U)$ then we have
$x\in\Ob(\Yg_U)$ if and only if $x|U_i\in\Ob(\Yg_{U_i})$ for all $i$.
\end{defn}

As an example, consider the stack $\Ello$ of generalised elliptic curves in the sense of 
\cite{delignerap}. Then $\Ell\subseteq\Ello$ is a substack: Since a generalised elliptic
curve is an elliptic curve if and only if it is smooth, condition i) of Definition \ref{def3} holds
because smoothness is stable under base change and condition ii) holds because smoothness if 
$fpqc$ local on the base.\\
\begin{defn}\label{def4}
A $1$-morphism $f:\Xg\longrightarrow\Yg$ of stacks is an epimorphism if for every $U\in\Ob(\Aff)$
and $y\in\Ob(\Yg_U)$ there exist a covering $\{ U_i\longrightarrow U\}$ in $\Aff$ and $x_i\in\Ob(\Xg_{U_i})$ such that $f_{U_i}(x_i)\simeq y|U_i$ for all $i$.
\end{defn}

A $1$-morphism of stacks is an isomorphism if and only if it is both a monomorphism and
an epimorphism \cite{CA}, Corollaire 3.7.1. This fact can also be understood from a homotopy theoretic
point of view \cite{hollander}, Corollary 8.16.\\
A fundamental insight is that many of the methods of algebraic geometry can be generalised to
apply to a suitable class of stacks. In order to define this class, we first have to explain 
the concept of representable $1$-morphisms of stacks which in turn needs the notion of 
algebraic spaces:\\
Algebraic spaces are a generalisation of schemes. The reader unfamiliar with them can, for the 
purpose of reading this paper, safely replace algebraic spaces by schemes throughout. We have
to mention them anyway in order to confirm with our main reference \cite{CA}. Algebraic spaces 
were invented by M. Artin and we decided not to try to give any short account of the main 
ideas underlying this master piece of algebraic geometry but rather refer the reader to
\cite{artin69} for an introduction and to \cite{knutson} as the standard technical reference.\\
We can now proceed on our way towards defining algebraic stacks.

\begin{defn}\label{def5}
A $1$-morphism $f:\Xg\longrightarrow\Yg$ of stacks is representable if for any $U\in\Aff$ with 
a $1$-morphism $U\longrightarrow\Yg$ the fibre product $\Xg\times_{\Yg} U$ is an algebraic space.
\end{defn}

Here, we refer the reader to \cite{CA}, 3.3 for the notion of finite limit for stacks.\\
Now let $P$ be a suitable property of morphisms of algebraic spaces, e.g. being an open or closed
immersion, being affine or being (faithfully) flat, see \cite{CA}, 3.10 for a more
exhaustive list. We say that a representable $1$-morphism $f:\Xg\longrightarrow\Yg$ of stacks
has the property $P$ if for every $U\in\Aff$ with a $1$-morphism $g:U\longrightarrow\Yg$, forming the cartesian
diagram
\[ \xymatrix{ \Xg\ar[r]^f& \Yg\\ \Xg\times_{\Yg} U\ar[u]\ar[r]^(.6){f'}& U \ar[u] ,} \]
the resulting morphism $f'$ between algebraic spaces has the property $P$.\\
As an example, let us check that the inclusion $\Ell\subseteq\Ello$ is an open immersion: To give
$U\in\Aff$ and a morphism $U\longrightarrow \Ello$ is the same as to give a generalised elliptic
curve $\pi:E\longrightarrow U$. Then $\Ell\times_{\Ello} U\longrightarrow U$ is the inclusion
of the complement of the image under $\pi$ of the non-smooth locus of $\pi$ and hence is an open
subscheme of $U$.\\
\begin{defn}\label{def6}
A stack $\Xg$ is algebraic if the diagonal $1$-morphism $\Xg\longrightarrow\Xg\times\Xg$ is representable
and affine and there is an affine scheme $U$ and a faithfully flat $1$-morphism $P: U\longrightarrow\Xg$.
\end{defn}

See section \ref{rigidstax} for further discussion.\\
A convenient way of constructing stacks is by means of groupoid objects. Let $(X_0,X_1)$ 
be a groupoid object in $\Aff$, i.e. a Hopf algebroid, see section \ref{stacksandhopf}.
Then $(X_0,X_1)$ determines a presheaf of groupoids on $\Aff$ and the corresponding category
fibred in groupoids $\Xg'$ is easily seen to satisfy condition i) of Definition \ref{def2}
for being a stack but not, in general, condition ii). There is a canonical way to pass from
$\Xg'$ to a stack $\Xg$ \cite{CA}, Lemme 3.2 which can also be interpreted as a fibrant 
replacement in a suitable model structure on presheaves of groupoids \cite{hollander}.\\
We provisionally define the stack of formal groups $\Xg_{FG}$ to be the stack associated with
the Hopf algebroid $(\MU_*,\MU_*\MU[u^{\pm 1}])$. Then $\Xg_{FG,U}'$ is the groupoid of 
formal group laws over $U$ and their (not necessarily strict) isomorphisms. A priori, it is 
unclear what the fibre categories $\Xg_{FG,U}$ are and in fact we will have to proceed differently
in section \ref{stackformalgroups}: We first construct a stack $\Xg_{FG}$ directly and then prove
that it is the stack associated with $(\MU_*,\MU_*\MU[u^{\pm 1}])$. \\
Note that there is a canonical $1$-morphism $\Spec(\MU_*)\longrightarrow\Xg_{FG}$. The following is a special case of Proposition \ref{landweberflat}.
\begin{prop} A $\MU_*$-algebra $R$ is Landweber exact if and only if the composition
$\Spec(R)\longrightarrow\Spec(\MU_*)\longrightarrow\Xg_{FG}$ is flat.
\end{prop}

\section{Algebraic stacks and flat Hopf algebroids}\label{stacksandhopf}
In this section we explain the relation between
flat Hopf algebroids and their categories of comodules and a certain class
of stacks and their categories of quasi-coherent sheaves of modules.

\subsection{The $2$-category of flat Hopf algebroids}\label{flatha}

We refer to \cite{R1}, Appendix A for the notion of a (flat) Hopf algebroid. To give
a Hopf algebroid $(A,\Gamma)$ is equivalent to giving $(X_0:=\Spec(A),
X_1:=\Spec(\Gamma))$ as a groupoid in affine schemes \cite{CA}, 2.4.3 and 
we will formulate most results involving Hopf algebroids this way.\\
Recall that this means that $X_0$ and $X_1$ are affine schemes and that we are given
morphisms $s,t:X_1\longrightarrow X_0$ (source and target), $\epsilon:X_0
\longrightarrow X_1$ (identity), $\delta:\cart{X_1}{s,X_0,t}{X_1}\longrightarrow X_1$ (composition) and $i:X_1\longrightarrow X_1$ (inverse) verifying suitable identities. The corresponding maps of rings are denoted $\eta_L,\eta_R$ (left-
and right unit), $\epsilon$ (augmentation), $\Delta$ (comultiplication) and
$c$ (antipode).\\
The $2$-category of flat Hopf algebroids ${\cal H}$ is defined as follows.
Objects are Hopf algebroids $(X_0,X_1)$ such that $s$ and $t$ are flat (and thus faithfully flat because they allow $\epsilon$ as a right inverse).
A $1$-morphism of flat Hopf algebroids from $(X_0,X_1)$ to $(Y_0,Y_1)$
is a pair of morphisms of affine schemes $f_i:X_i\longrightarrow Y_i$ ($i=0,1$) commuting with
all the structure. The composition of $1$-morphisms is component wise.
Given two $1$-morphisms $(f_0,f_1),(g_0,g_1):(X_0,X_1)\longrightarrow (Y_0,Y_1)$, 
a $2$-morphism $c:(f_0,f_1)\longrightarrow (g_0,g_1)$ is a morphism of affine schemes
$c:X_0\longrightarrow Y_1$ such that $sc=f_0, tc=g_0$ and the diagram
\[ \xymatrix{ X_1\ar[r]^-{(g_1,cs)}\ar[d]_{(ct,f_1)} & Y_1{\times\atop s,Y_0,t}Y_1\ar[d]^{\delta} \\
Y_1{\times\atop s,Y_0,t}Y_1\ar[r]^-{\delta} & Y_1} \]
commutes. For $(f_0,f_1)=(g_0,g_1)$ the identity $2$-morphism is given by
$c:=\epsilon f_0$. Given two $2$-morphisms $\xymatrix{(f_0,f_1)
\ar[r]^c & (g_0,g_1)\ar[r]^{c^{'}} & (h_0,h_1)}$ their composition is defined as
\[ c^{'}\circ c:\xymatrix{X_0\ar[r]^-{(c^{'},c)} & Y_1{\times\atop s,Y_0,t}Y_1 \ar[r]^-{\delta} & Y_1}. \]
One checks that the above definitions make ${\cal H}$ a $2$-category which is 
in fact clear because, except for the flatness of $s$ and $t$, they are 
merely a functorial way of stating the axioms of a groupoid, a functor and 
a natural transformation. For technical reasons we will sometimes consider Hopf
algebroids for which $s$ and $t$ are not flat.\\

\subsection{The $2$-category of rigidified algebraic stacks}\label{rigidstax}

From Definition \ref{def2} one sees that any $1$-morphism of algebraic stacks from an algebraic space to an 
algebraic stack is representable and affine, c.f. the proof of \cite{CA}, Corollaire 3.13.
In particular, the condition in Definition \ref{def2} that $P$ be
faithfully flat makes sense. 
By definition, every algebraic stack is quasi-compact, hence so is any
$1$-morphism between algebraic stacks \cite{CA}, D\'efinition 4.16, Remarques 4.17.
One can check that finite limits and colimits of algebraic stacks 
are again algebraic stacks. If $\Ug\stackrel{i}{\hookrightarrow}\Xg$ is a quasi-compact open immersion of stacks and $\Xg$ is algebraic then the stack $\Ug$
is algebraic as one easily checks. In general, an open substack of an algebraic stack need not be algebraic, see the introduction of section \ref{ringchange}.\\
A morphism $P$ as in Definition \ref{def2} is called a presentation of $\Xg$. As far as we are aware,
the above definition of ``algebraic'' is due to P. Goerss \cite{G} and is
certainly motivated by the equivalence given in subsection \ref{equivalence}
below. We point out that the notion of ``algebraic stack'' well-establish 
in algebraic geometry \cite{CA}, D\'efinition 4.1 is different from the above. For example,
the stack associated with $(\BP_*,\BP_*\BP)$ in section \ref{ringchange} 
is algebraic in the above sense but not in the sense of algebraic geometry because its diagonal is not of finite type \cite{CA} Lemme 4.2.
Of course, in the following we will use the term ``algebraic stack'' in the sense defined above.\\
The $2$-category ${\cal S}$ of rigidified algebraic stacks is defined as follows.
Objects are presentations $P: X_0\longrightarrow \Xg$ as in Definition \ref{def2}. A $1$-morphism
from $P: X_0\longrightarrow \Xg$ to $Q: Y_0\longrightarrow \Yg$ is a pair
consisting of $f_0:X_0\longrightarrow Y_0$ in $\Aff$ and a $1$-morphism of stacks $f:\Xg\longrightarrow \Yg$ such that the diagram

\[ \xymatrix{ X_0\ar[r]^{f_0}\ar[d]_P & Y_0\ar[d]^Q \\
\Xg \ar[r]_f &\Yg }\]

is $2$-commutative. The composition of $1$-morphisms is component wise. Given 
$1$-morphisms $(f_0,f), (g_0,g):(X_0\longrightarrow \Xg)\longrightarrow
(Y_0\longrightarrow \Yg)$ a $2$-morphism in ${\cal S}$ from $(f_0,f)$ to
$(g_0,g)$ is by definition a $2$-morphism from $f$ to $g$ in the $2$-category
of stacks \cite{CA}, 3.

\hspace{.5 cm}

\subsection{The equivalence of ${\cal H}$ and ${\cal S}$}\label{equivalence}

We now establish an equivalence of $2$-categories between ${\cal H}$ and 
${\cal S}$. We define a functor
$K:{\cal S}\longrightarrow {\cal H}$ as follows.\\
\[ K(\xymatrix{X_0\ar[r]^P & \Xg}):=(X_0, X_1:=\cart{X_0}{P,\Xg,P}{X_0}) \]
has a canonical structure of groupoid \cite{CA}, Proposition 3.8, $X_1$ is affine because $X_0$ is affine and 
$P$ is representable and affine and the projections $s,t:\xymatrix{X_1\ar@<.5ex>[r]\ar@<-.5ex>[r] & X_0}$ are flat
because $P$ is. Thus $(X_0,X_1)$ is a flat Hopf algebroid. If $(f_0,f): \xymatrix{(X_0\ar[r]^P & \Xg)}\longrightarrow
\xymatrix{(Y_0\ar[r]^Q & \Yg)}$ is a $1$-morphism in ${\cal S}$
we define $K((f_0,f)):=(f_0,f_0\times f_0)$. If we have $1$-morphisms
$(f_0,f),(g_0,g): \xymatrix{(X_0\ar[r]^P & \Xg)}\longrightarrow
\xymatrix{(Y_0\ar[r]^Q & \Yg)}$ in ${\cal S}$ and a $2$-morphism 
$(f_0,f)\longrightarrow (g_0,g)$ then we have by definition a $2$-morphism
$\xymatrix{f \ar[r]^{\Theta} & g}: \Xg\longrightarrow \Yg$. In particular,
we have $\Theta_{X_0}:\Ob(\Xg_{X_0})\longrightarrow\Mor(\Yg_{X_0})=\Hom_{\Aff}(X_0,Y_1)$ and we define $K(\Theta):=\Theta_{X_0}(\id_{X_0})$.
One checks that $K:{\cal S}\longrightarrow {\cal H}$ is a $2$-functor.\\
We define a $2$-functor $G:{\cal H}\longrightarrow {\cal S}$ as follows.
On objects we put $G((X_0,X_1)):=(X_0\stackrel{can}{\longrightarrow} \Xg:=[ \xymatrix{ X_1\ar@<.5ex>[r]\ar@<-.5ex>[r] & X_0} ] )$, the stack associated
with the groupoid $(X_0,X_1)$ together with its canonical presentation \cite{CA}, 3.4.3; identify the $X_i$ with the flat sheaves they represent to consider them as ``S-espaces'', see also subsection \ref{epimonic}. Then $G((X_0,X_1))$ is 
a rigidified algebraic stack: Saying that the diagonal of $\Xg$ is
representable and affine means that for any algebraic space $X$ and morphisms $x_1,x_2:X\longrightarrow\Xg$ the sheaf $\Isom_X(x_1,x_2)$ on $X$ is representable by an affine $X$-scheme. This problem is local in the {\em fpqc} topology on $X$ because affine morphisms satisfy effective descent in the {\em fpqc} topology \cite{SGA1}, expos\'e VIII, Th\'eor\`eme 2.1. So we can assume that the $x_i$ lift to $X_0$ and the assertion follows because $(s,t): X_1\longrightarrow\cart{X_0}{S}{X_0}$ is affine. A similar argument shows that $P:X_0\longrightarrow\Xg$ is representable and faithfully flat because $s$ and $t$ are faithfully flat.\\
Given a $1$-morphism $(f_0,f_1):(X_0,X_1)\longrightarrow (Y_0,Y_1)$
in ${\cal H}$ there is a unique $1$-morphism $f:\Xg\longrightarrow\Yg$ making
\[ \xymatrix{ X_1\ar@<-.5ex>[r]\ar@<.5ex>[r]\ar[d]_{f_1} & X_0 \ar[r]^-P\ar[d]_{f_0} & \Xg \ar@{.>}[d]^f \\
Y_1\ar@<-.5ex>[r]\ar@<.5ex>[r] & Y_0 \ar[r]^-Q & \Yg } \]

$2$-commutative \cite{CA}, proof of Proposition 4.18 and we define $G((f_0,f_1)):=f$.\\
Given a $2$-morphism $c:X_0\longrightarrow Y_1$ from the $1$-morphism $(f_0,f_1):(X_0,X_1)\longrightarrow (Y_0,Y_1)$ to the $1$-morphism $(g_0,g_1):(X_0,X_1)\longrightarrow (Y_0,Y_1)$ in ${\cal H}$ we have a diagram 
\[ \xymatrix{ X_1\ar@<-.5ex>[r]\ar@<.5ex>[r]\ar@<-.5ex>[d]_{f_1}\ar@<.5ex>[d]^{g_1} & X_0 \ar[r]^-P\ar@<-.5ex>[d]_{f_0}\ar@<.5ex>[d]^{g_0}  & \Xg \ar@/^/[d]^-g\ar@/_/[d]_-f \\
Y_1\ar@<-.5ex>[r]\ar@<.5ex>[r] & Y_0 \ar[r]^-Q & \Yg } \]
and need to construct a $2$-morphism $\Theta=G(c):f\longrightarrow g$ in the 
$2$-category of stacks. We will do this in some detail because we omit numerous similar
arguments.\\
Fix $U\in\Aff$, $x\in\Ob(\Xg_U)$ and a representation of $x$ as in \cite{CA}, proof of Lemme 3.2 \[ (U^{'}\longrightarrow U,\, x^{'}:U^{'}\longrightarrow X_0,\, U^{''}:=\xymatrix{ U^{'} \ar@{}[r]|{\times}_-U & U^{'}}\stackrel{\sigma}{\longrightarrow} X_1), \]
i.e. $U^{'}\longrightarrow U$ is a cover in $\Aff$, $x^{'}\in X_0(U^{'})=\Hom_{\Aff}(U^{'},
X_0)$ and $\sigma$ is a descent datum for $x^{'}$ with respect to the cover 
$U^{'}\longrightarrow U$. Hence, denoting by $\pi_1,\pi_2:U^{''}\longrightarrow U^{'}$
and $\pi: U^{'}\longrightarrow U$ the projections, we have $\sigma:\pi_1^*x^{'}\stackrel{\sim}{\longrightarrow}\pi_2^*x^{'}$ in $\Xg_{U^{''}}$, i.e. $x^{'}\pi_1=s\sigma$ and
$x^{'}\pi_2=t\sigma$. Furthermore, $\sigma$ satisfies a cocycle condition which we do not 
spell out.\\
We have to construct a morphism \[ \Theta_x:f(x)\longrightarrow g(x) \mbox{ in }\Yg_U\]
which we do by descent from $U^{'}$ as follows. We have a morphism \[ \pi^*(f(x))
=f(\pi^*(x)=x^{'})=f_0x^{'}\stackrel{\phi^{'}}{\longrightarrow}\pi^*(g(x))=g_0x^{'}
\mbox{ in }\Yg_{U^{'}}\]
given by $\phi^{'}:=cx^{'}:U^{'}\longrightarrow Y_1$. We also have a diagram 
\[ \xymatrix{ \pi_1^*(\pi^*(f(x)))=f_0x^{'}\pi_1 \ar[r]^-{\pi_1^*(\phi^{'})} 
\ar[d]_-{\sigma_f} & \pi_1^*(\pi^*(g(x)))=g_0x^{'}\pi_1 \ar[d]^{\sigma_g} \\
\pi_2^*(\pi^*(f(x)))=f_0x^{'}\pi_2 \ar[r]^-{\pi_2^*(\phi^{'})} & \pi_2^*(\pi^*(g(x)))=g_0x^{'}\pi_2} \]

in $\Yg_{U^{''}}$ where $\sigma_f$ and $\sigma_g$ are descent isomorphisms for $f(x^{'})$ and $g(x^{'})$ given by 
$\sigma_f=f_1\sigma$ and $\sigma_g=g_1\sigma$. We check that this diagram commutes 
by computing in $\Mor(\Yg_{U^{''}})$: \[ \sigma_g\circ \pi_1^*(\phi^{'})=
\delta_Y(g_1\sigma,cx^{'}\pi_1)=\delta_Y(g_1\sigma,cs\sigma)=\delta_Y(g_1,cs)\sigma
\stackrel{(*)}{=}\]
\[ =\delta_Y(ct,f_1)\sigma=\delta_Y(ct\sigma,f_1\sigma)=\delta_Y(cx^{'}\pi_2,f_1\sigma)=\pi_2^*(\phi^{'})\circ\sigma_f.\]
Here $\delta_Y$ is the composition of $(Y_0,Y_1)$ and in $(*)$ we used the commutative 
square in the definition of $2$-morphisms in ${\cal H}$.\\
So $\phi^{'}$ is compatible with descent data and thus descents to the desired $\Theta_x:
f(x)\longrightarrow g(x)$. We omit the verification that $\Theta_x$ is independent of
the chosen representation of $x$ and natural in $x$ and $U$.
One checks that $G:{\cal H}\longrightarrow{\cal S}$ is a $2$-functor.\\

\begin{theorem}\label{equivofSandH}
The above $2$-functors $K:{\cal S}\longrightarrow {\cal H}$ and $G:{\cal H}\longrightarrow {\cal S}$ are inverse equivalences.
\end{theorem}

\begin{proof} We have $G\circ K(X_0\stackrel{P}{\longrightarrow}\Xg)=(\xymatrix{ X_0 \ar[r]^-{can}
 & [ X_0{\times\atop\Xg}X_0 \ar@<-.5ex>[r]\ar@<.5ex>[r] & X_0 ]})$ and there 
is a unique $1$-isomorphism $\nu_P:\xymatrix{[ X_0{\times\atop\Xg}X_0 \ar@<-.5ex>[r]\ar@<.5ex>[r] & X_0 ]}\longrightarrow\Xg$ with
$\nu_p\circ can=P$ \cite{CA}, Proposition 3.8. One checks that this defines an isomorphism of $2$-functors $G\circ K\stackrel{\simeq}{\longrightarrow}\id_{\cal S}$.\\
Next we have $K\circ G(X_0,X_1)=(X_0,\cart{X_0}{P,\Xg,P}{X_0})$, where 
$(X_0\stackrel{P}{\longrightarrow}\Xg)=G(X_0,X_1)$, and $X_1\simeq\cart{X_0}{P,\Xg,P}{X_0}$ \cite{CA}, 3.4.3 and one checks that this defines an isomorphism of $2$-functors $\id_{\cal H}\stackrel{\simeq}{\longrightarrow} K\circ G$.
\end{proof}
In the following, given a flat Hopf algebroid $(X_0,X_1)$, we will refer
to $G((X_0,X_1))$ simply as the (rigidified) algebraic stack associated with 
$(X_0,X_1)$.\\
The forgetful functor from rigidified algebraic stacks to algebraic stacks 
is not full but we have the following.

\begin{prop}\label{chain} If $(X_0,X_1)$ and $(Y_0,Y_1)$ are flat Hopf algebroids
with associated rigidified algebraic stacks $P:X_0\longrightarrow\Xg$ and 
$Q:Y_0\longrightarrow\Yg$ and $\Xg$ and $\Yg$ are $1$-isomorphic as stacks then
there is a chain of $1$-morphisms of flat Hopf algebroids from $(X_0,X_1)$
to $(Y_0,Y_1)$ such that every morphism in this chain induces a $1$-isomorphism
on the associated algebraic stacks.
\end{prop}
\begin{remark} This result implies Theorem 6.5 of \cite{H1}: By Theorem \ref{image} below, the assumptions of {\em loc. cit.} imply that the flat 
Hopf algebroids $(B,\Gamma_B)$ and $(B^{'},\Gamma_{B^{'}})$ considered there
have the same open substack of the stack of formal groups as their associated stack. So they are connected by a chain of weak equivalences by Proposition
\ref{chain}, see Remark \ref{tannakabla} for the notion of weak equivalence.
\end{remark}

\begin{proof} Let $f:\Xg\longrightarrow\Yg$ be a $1$-isomorphism of stacks and form the cartesian diagram
\[ \xymatrix{ X_1^{'}\ar[r]_-{f_1} \ar@<-.5ex>[d] \ar@<.5ex>[d] & Y_1 \ar@<-.5ex>[d] \ar@<.5ex>[d] \\
X_0^{'} \ar[r]_-{f_0} \ar[d]_-{P^{'}} & Y_0 \ar[d]^-{Q} \\
\Xg \ar[r]_-{f} & \Yg. } \]
To be precise, the upper square is cartesian for either both source or both target morphisms.
Then $(f_0,f_1)$ is a $1$-isomorphism of flat Hopf algebroids.
Next, $Z:=\cart{X_0^{'}}{P^{'},\Xg,P}{X_0}$ is an affine scheme because 
$X_0^{'}$ is and $P$ is representable and affine. The obvious $1$-morphism
$Z\longrightarrow\Xg$ is representable, affine and faithfully flat because
$P$ and $P^{'}$ are. Writing $W:=\cart{Z}{\Xg}{Z}\simeq\cart{X_1^{'}}{\Xg}{X_1}$
we have that $\Xg\simeq [ \xymatrix{ W \ar@<.5ex>[r]\ar@<-.5ex>[r] & Z} ]$
by the flat version of \cite{CA}, Proposition 4.3.2. There are obvious $1$-morphisms of
flat Hopf algebroids $(Z,W)\longrightarrow (X_0^{'},X_1^{'})$ and 
$(Z,W)\longrightarrow (X_0,X_1)$ covering $\id_{\Xg}$ (in particular inducing
an isomorphism on stacks) and we get the sought for chain as $(Y_0,Y_1)\longleftarrow (X_0^{'},X_1^{'})\longleftarrow (Z,W)\longrightarrow (X_0,X_1)$.
\end{proof}

\vspace{1.5cm}

\subsection{Comodules and quasi-coherent sheaves of modules}\label{modules}

Let $(A,\Gamma)$ be a flat Hopf algebroid with associated rigidified
algebraic stack $X_0=\Spec(A)\longrightarrow\Xg$. From Theorem \ref{equivofSandH} one would certainly expect that the category of $\Gamma$-comodules
has a description in terms of $X_0\longrightarrow\Xg$. In this subsection we
prove the key observation that this category does in fact only depend on $\Xg$
and not on the particular presentation $X_0\longrightarrow\Xg$, c.f. (\ref{boing}) below.\\
For basic results concerning the category $\Modq{\Xg}$ of quasi-coherent
sheaves of modules on an algebraic stack $\Xg$ we refer the reader to \cite{CA}, 13.\\
Fix a rigidified algebraic stack $X_0\stackrel{P}{\longrightarrow}\Xg$ corresponding 
by Theorem \ref{equivofSandH} to the flat Hopf algebroid $(X_0=\Spec(A),X_1=\Spec(\Gamma))$ with structure morphisms $s,t:X_1\longrightarrow X_0$. As $P$ is affine it is in particular quasi-compact, hence {\em fpqc}, and thus of effective cohomological descent for quasi-coherent
modules \cite{CA}, Th\'eor\`eme 13.5.5,i). In particular, $P^*$ induces an equivalence
\[ P^*:\ModX\stackrel{\simeq}{\longrightarrow}\{F\in\ModXz+\mbox{ descent data}\},\]
c.f. \cite{BLR}, Chapter 6 for similar examples of descent.  A descent datum on $F\in\ModXz$ is an isomorphism
$\alpha:s^*F\longrightarrow t^*F$ in $\ModXo$ satisfying a cocycle condition. Giving
$\alpha$ is equivalent to giving either its adjoint $\psi_l:F\longrightarrow s_*t^*F$
or the adjoint of $\alpha^{-1}$, $\psi_r: F\longrightarrow t_*s^*F$. Writing 
$M$ for the $A$-module corresponding to $F$, $\alpha$ corresponds to an isomorphism 
$\ocart{\Gamma}{\eta_L,A}{M}\longrightarrow\ocart{\Gamma}{\eta_R,A}{M}$ of 
$\Gamma$-modules and $\psi_r$ and $\psi_l$ correspond respectively to morphisms
$M\longrightarrow \Gamma\otimes_{\eta_R,A} M$ and $M\longrightarrow M\otimes_{A,\eta_L}\Gamma$
 of $A$-modules.
One checks that this is a $1$-$1$ correspondence between descent data on $F$ 
and left- (respectively right-)$\Gamma$-comodule structures on $M$. For example,
the cocycle condition for $\alpha$ corresponds to the coassociativity of the coaction.
In the following we will work with left-$\Gamma$-comodules exclusively and simply call them 
$\Gamma$-comodules. The above construction then provides an explicit equivalence 
\begin{equation}\label{boing}
\ModX\stackrel{\simeq}{\longrightarrow} \Gamma\mbox{-comodules.}
\end{equation}
This can also be proved using the Barr-Beck theorem, \cite{P}, 3.22.\\
The identification of $\ModX$ with $\Gamma$-comodules allows to (re)understand a number of results 
on $\Gamma$-comodules from the stack theoretic point of view and we now give a short list of
such applications which we will use later.\\
The adjunction $(P^*,P_*):\ModX\longrightarrow\ModXz$ corresponds to the 
forgetful functor from $\Gamma$-comodules to $A$-modules, respectively to the functor
``induced/extended comodule''. The structure sheaf ${\cal O}_{\Xg}$ corresponds to the trivial 
$\Gamma$-comodule $A$, hence taking the primitives of a $\Gamma$-comodule 
(i.e. the functor $\Hom_{\Gamma}(A,\cdot)$ from $\Gamma$-comodules to abelian groups) corresponds to $\Hom_{{\cal O}_\Xg}({\cal O}_\Xg,\cdot)=H^0(\Xg,\cdot)$ and thus $\Ext^*_{\Gamma}(A,\cdot)$ corresponds to quasi-coherent
cohomology $H^*(\Xg,\cdot)$. Another application of (\ref{boing})
is the following correspondence between closed substacks and invariant ideals:\\
By \cite{CA}, Application 14.2.7 there is a $1$-$1$ correspondence between closed substacks $\Zg\subseteq\Xg$ and quasi-coherent ideal sheaves $\Ih\subseteq\Oh_{\Xg}$ under which  $\Oh_{\Zg}\simeq
\Oh_{\Xg}/\Ih$ and by (\ref{boing}) these $\Ih$ correspond to $\Gamma$-subcomodules $I\subseteq A$, i.e.
invariant ideals. In this situation, the diagram

\[ \xymatrix{ \Spec(\Gamma/I\Gamma)\ar[r]\ar@<-.5ex>[d]\ar@<.5ex>[d] & \Spec(\Gamma)\ar@<-.5ex>[d]\ar@<.5ex>[d] \\
\Spec(A/I) \ar[r]\ar[d] & \Spec(A)\ar[d] \\ \Zg\ar[r] & \Xg }\]

is cartesian. Note that the Hopf algebroid $(A/I,\Gamma/I\Gamma)$ is induced from $(A,\Gamma)$
by the map $A\longrightarrow A/I$ because $A/I\otimes_A\Gamma\otimes_A A/I\simeq
\Gamma/(\eta_LI+\eta_RI)\Gamma=\Gamma/I\Gamma$ since $I$ is invariant.\\
We conclude this subsection by giving a finiteness result for quasi-coherent
sheaves of modules. Let $\Xg$ be an algebraic stack. We say that $\Fh\in\Modq{\Xg}$
if {\em finitely generated} if there is a presentation $P:X_0=\Spec(A)\longrightarrow\Xg$
such that the $A$-module corresponding to $P^*\Fh$ is finitely generated. If
$\Fh$ is finitely generated then for any presentation $P:X_0^{'}=\Spec(A^{'})\longrightarrow\Xg$ the $A^{'}$-module corresponding to $P^{'*}\Fh$ is finitely generated as one sees using \cite{Bou}, I, \S3, Proposition 11. 
\begin{prop}\label{fgenha} Let $(A,\Gamma)$ be a flat Hopf algebroid,
$M$ a $\Gamma$-comodule and $M^{'}\subseteq M$ a finitely generated $A$-submodule. Then $M^{'}$ is contained in a $\Gamma$-subcomodule of $M$ which is finitely 
generated as an $A$-module.
\end{prop}
\begin{proof} \cite{Thorsten}, Proposition 5.7. \end{proof}
Note that in this result, ``finitely generated'' cannot be strengthened
to ``coherent'' as is shown by the example of the simple $\BP_*\BP$-comodule
$\BP_*/(v_0,\ldots)$ which is not coherent as a $\BP_*$-module.
\begin{prop}\label{fgen} Let $\Xg$ be an algebraic stack. Then any 
$\Fh\in\Modq{\Xg}$ is the filtering union of its finitely generated quasi-coherent subsheaves.
\end{prop}
\begin{proof} Choose a presentation of $\Xg$ and apply Proposition \ref{fgenha}
to the resulting flat Hopf algebroid. \end{proof}

This result may be compared with \cite{CA}, Proposition 15.4.

\section{Tannakian results}\label{morphisms}

In \cite{Lurie}, J. Lurie considers a Tannakian correspondence for "geometric" stacks which
are exactly those stacks that are algebraic {\em both} in the sense of \cite{CA}, D\'efinition 
4.1 and in the sense of Definition \ref{def6}. He shows that associating to such
a stack $\Xg$ the category $\ModX$ is a fully faithful $2$-functor. The recognition problem, 
i.e. giving an intrinsic characterisation of the categories $\ModX$, remains open but see
\cite{tannaka} for a special case.\\
The usefulness of a Tannakian correspondence stems from being able to relate notions of
linear algebra, pertaining to the categories $\ModX$ and their morphisms, to geometric
notions, pertaining to the stacks and their morphisms. See \cite{DMOS},Propositions 2.20-29 for examples of
this in the special case that $\Xg=BG$ is the classifying stack of a linear algebraic group $G$.
This relation can be studied without having solved the recognition problem and we do so in
the present section, i.e. we relate properties of $1$-morphisms $(f_0,f_1)$ 
of flat Hopf algebroids to properties of the induced morphism $f:\Xg
\longrightarrow\Yg$ of algebraic stacks and the adjoint pair
$(f^*,f_*):\Modq{\Xg}\longrightarrow\Modq{\Yg}$ of functors.

\subsection{The epi/monic factorisation}\label{epimonic}

Every $1$-morphism of stacks factors canonically into an epimorphism
followed by a monomorphism and in this subsection we explain the analogous 
result for (flat) Hopf algebroids. In particular, this will explain
the stack theoretic meaning of the construction of an induced Hopf
algebroid, c.f. \cite{H1}, beginning of section 2.\\
By a {\em flat sheaf} we will mean a set valued sheaf on the site $\Aff$.
The topology of $\Aff$ is subcanonical, i.e. every representable 
presheaf is a sheaf. We can thus identify the category underlying $\Aff$ 
with a full subcategory of the category of flat sheaves.\\
Every $1$-morphism $f:\Xg\longrightarrow\Yg$ of stacks factors canonically
$\Xg\longrightarrow\Xg^{'}\longrightarrow\Yg$ into an epimorphism followed
by a monomorphism \cite{CA}, Proposition 3.7. The stack $\Xg^{'}$ is determined up to
unique $1$-isomorphism and is called the image of $f$.\\
For a $1$-morphism $(f_0,f_1):(X_0,X_1)\longrightarrow (Y_0,Y_1)$ 
of flat Hopf algebroids we introduce
\begin{eqnarray}\label{ab}
 & \alpha:=t\pi_2 : \cart{X_0}{f_0,Y_0,s}{Y_1}\longrightarrow Y_0\mbox{ and } & \\
 &  \beta:=(s,f_1,t): X_1\longrightarrow\cartt{X_0}{f_0,Y_0,s}{Y_1}{t,Y_0,f_0}{X_0}. &\nonumber  
\end{eqnarray}
The $1$-morphism $f:\Xg\longrightarrow \Yg$ induced by $(f_0,f_1)$ on algebraic stacks is an epimorphism if and only if
$\alpha$ is an epimorphism of flat sheaves as is clear from Definition \ref{def4}. On the other hand, $f$ is a monomorphism if and only
if $\beta$ is an isomorphism, as is easily checked.\\
Writing $X_1^{'}:=\cartt{X_0}{f_0,Y_0,s}{Y_1}{t,Y_0,f_0}{X_0}$, $(f_0,f_1)$ factors as
\[ \xymatrix{ X_1\ar[r]^{f_1^{'}:=\beta}\ar@<-.5ex>[d]\ar@<.5ex>[d]
 & X_1^{'}\ar[r]^{\pi_2}\ar@<-.5ex>[d]_-{\pi_1}\ar@<.5ex>[d]^-{\pi_3} & Y_1\ar@<-.5ex>[d]\ar@<.5ex>[d] \\
X_0 \ar[r]^{f_0^{'}:=\id_{X_0}} & X_0 \ar[r]^{f_0} & Y_0 } \]

and the factorisation of $f$ induced by this is the epi/monic factorisation.
Note that even if $(X_0,X_1)$ and $(Y_0,Y_1)$ are flat Hopf algebroids, $(X_0,X_1^{'})$ does not have to be flat.\\
We refer to $(X_0,X_1^{'})$ as the Hopf algebroid induced from $(Y_0,Y_1)$ 
by $f_0$.

\subsection{Flatness and isomorphisms}\label{flatandiso}

The proof of the next result will be given at the end of this 
subsection. The equivalence of $ii)$ and $iii)$ is equivalent to Theorem 
6.2 of \cite{H1} but we will obtain refinements of it below, see Proposition
\ref{epiofstack} and Proposition \ref{monothm}.

\begin{theorem}\label{isothm} Let $(f_0,f_1):(X_0,X_1)\longrightarrow (Y_0,Y_1)$
be a $1$-morphism of flat Hopf algebroids with associated morphisms $\alpha$
and $\beta$ as in (\ref{ab}) and inducing $f:\Xg\longrightarrow\Yg$ on algebraic stacks.
Then the following are equivalent:\\
i) $f$ is a $1$-isomorphism of stacks.\\
ii) $f^*:\Modq{\Xg}\longrightarrow\Modq{\Yg}$ is an equivalence.\\
iii) $\alpha$ is faithfully flat and $\beta$ is an isomorphism.
\end{theorem}

\begin{remark}\label{tannakabla} 
This result shows that weak equivalences as defined in \cite{H2}, Definition 1.1.4
are exactly those $1$-morphisms of flat Hopf algebroids which induce 
$1$-isomorphisms on the associated algebraic stacks.\\
\end{remark}

We next give two results about the flatness of morphisms.

\begin{prop}\label{flatness1} Let $(f_0,f_1):(X_0,X_1)\longrightarrow (Y_0,Y_1)$ be a $1$-morphism of flat Hopf algebroids, $P:X_0\longrightarrow\Xg$ and 
$Q:Y_0\longrightarrow\Yg$ the associated rigidified algebraic stacks and
$f:\Xg\longrightarrow\Yg$ the induced $1$-morphism of algebraic stacks. Then the following 
are equivalent:\\
i) $f$ is (faithfully) flat.\\ 
ii) $f^*:\Modq{\Yg}\longrightarrow\Modq{\Xg}$ is exact (and faithful).\\
iii) $\alpha:=t\pi_2 :\cart{X_0}{f_0,Y_0,s}{Y_1}\longrightarrow Y_0$
is (faithfully) flat.\\
iv) The composition $X_0\stackrel{P}{\longrightarrow}\Xg\stackrel{f}{\longrightarrow}\Yg$ is (faithfully) flat.
\end{prop}

\begin{proof} The equivalence of i) and ii) holds by definition, the one of 
i) and iv) holds because $P$ is $fpqc$ and being (faithfully) flat is a local
property for the {\em fpqc} topology.
Abbreviating $Z:=\cart{X_0}{f_0,Y_0,s}{Y_1}$ we have a cartesian diagram 
\[ \xymatrix{ Z \ar[rr]^-{\alpha}\ar[d]_-{\pi_1} & &
Y_0 \ar[d]^-Q \\ 
X_0 \ar[r]^-P \ar@{.>}[rru]^{f_0} & \Xg \ar[r]^-f  & \Yg } \]
which, as $Q$ is $fpqc$, shows that iv) and iii) are equivalent. We check that
this diagram is in fact cartesian by computing:
\[ \cart{X_0}{fP,\Yg,Q}{Y_0}=\cart{X_0}{Qf_0,\Yg,Q}{Y_0}\simeq \]
\[ \simeq\cartt{X_0}{f_0,Y_0,\id}{Y_0}{Q,\Yg,Q}{Y_0}\simeq\cart{X_0}{f_0,Y_0,s}{Y_1}=Z,\]
and under this isomorphism the projection onto the second factor corresponds to $\alpha$.
\end{proof}

\begin{prop}\label{flatness2} Let $(Y_0,Y_1)$ be a flat Hopf algebroid, $f_0:
X_0\longrightarrow Y_0$ a morphism in $\Aff$ and $(f_0,f_1):(X_0,
X_1:=\cartt{X_0}{f_0,Y_0,s}{Y_1}{t,Y_0,f_0}{X_0})\longrightarrow (Y_0,Y_1)$
the canonical $1$-morphism of Hopf algebroids from the induced
Hopf algebroid and $Q:Y_0\longrightarrow\Yg$ the rigidified algebraic stack
associated with $(Y_0,Y_1)$. Then the following are equivalent:\\
i) The composition $X_0\stackrel{f_0}{\longrightarrow}Y_0\stackrel{Q}{\longrightarrow}\Yg$ is (faithfully) flat.\\
ii) $\alpha:=t\pi_2 :\cart{X_0}{f_0,Y_0,s}{Y_1}\longrightarrow Y_0$
is (faithfully) flat.\\
If either of this maps is flat, then $(X_0,X_1)$ is a {\em flat} Hopf
algebroid.
\end{prop}

The last assertion of this Proposition does not admit a converse: 
For $(Y_0,Y_1)=(\Spec(\BP_*),\Spec(\BP_*\BP))$ and $X_0:=\Spec(\BP_*/I_n)\longrightarrow
Y_0$, the induced Hopf algebroid is flat but $X_0\longrightarrow\Yg$ is
not, c.f. subsection \ref{landweber}.

\begin{proof} The proof of the equivalence of i) and ii) is the same as 
in Proposition \ref{flatness1},
using that $Q$ is $fpqc$. Again denoting $Z:=\cart{X_0}{f_0,Y_0,s}{Y_1}$ one
checks that the diagram
\[ \xymatrix{ Z \ar[r]^-{\alpha} & Y_0\\
X_1 \ar[u] \ar[r]^-t & X_0 \ar[u]_-{f_0}} \]
is cartesian which implies the final assertion of the proposition because flatness
is stable under base change.
\end{proof}

\begin{prop}\label{inducediso} 
Let $(Y_0,Y_1)$ be a flat Hopf algebroid, $f_0:
X_0\longrightarrow Y_0$ a morphism in $\Aff$ such that the composition
$X_0\stackrel{f_0}{\longrightarrow}Y_0\stackrel{Q}{\longrightarrow}\Yg$ is
faithfully flat, where $Q:Y_0\longrightarrow\Yg$ is the rigidified algebraic 
stack associated with $(Y_0,Y_1)$. Let $(f_0,f_1):(X_0,X_1)\longrightarrow
(Y_0,Y_1)$ be the canonical $1$-morphism with $(X_0,X_1)$ the Hopf algebroid
induced from $(Y_0,Y_1)$ by $f_0$. Then $(X_0,X_1)$ is a flat Hopf 
algebroid and $(f_0,f_1)$ induces a $1$-isomorphism on the associated
algebraic stacks.
\end{prop}

\begin{proof} The $1$-morphism $f$ induced on the associated
algebraic stacks is a monomorphism as explained in subsection \ref{epimonic}. Proposition \ref{flatness2}
shows that $(X_0,X_1)$ is a flat Hopf algebroid and that $\alpha$ is faithfully flat, hence an epimorphism of flat sheaves. Thus $f$ is an epimorphism
of stacks as noted in subsection \ref{epimonic} and, finally, $f$ is a $1$-
isomorphism by \cite{CA}, Corollaire 3.7.1.
\end{proof}

We now start to take the module categories into consideration.
Given $f:X\longrightarrow Y$ in $\Aff$ we have an adjunction $\psi_f:
\id_{\Modq{Y}}\longrightarrow f_*f^*$.
We recognise the epimorphisms of representable flat
sheaves as follows.

\begin{prop}\label{epiofrep}
Let $f:X\longrightarrow Y$ be a morphism in $\Aff$. Then the following 
are equivalent:\\
i) $f$ is an epimorphism of flat sheaves.\\
ii) There is some $\phi:Z\longrightarrow X$ in $\Aff$ such
that $f\phi$ is faithfully flat.\\
If i) and ii) hold, then $\psi_f$ is injective.\\
If $f$ is flat, the conditions i) and ii) are equivalent to $f$ being
faithfully flat.

\end{prop}

As an example of a morphism satisfying the conditions of Proposition \ref{epiofrep}
without being flat one may take the unique morphism $\Spec(\Z)\sqcup\Spec(\F_p)\longrightarrow\Spec(\Z)$.

\begin{proof}
That i) implies ii) is seen by lifting $\id_Y\in Y(Y)$ after a suitable
faithfully flat cover $Z\longrightarrow Y$ to some $\phi\in X(Z)$.\\
To see that ii) implies i), fix some $U\in\Aff$ and $u\in Y(U)$ and form 
the cartesian diagram 
\[ \xymatrix{ Z\ar[r]^-{\phi} & X\ar[r]^-f & Y \\
W\ar[u]_-v \ar[rr] & & U\ar[u]_-u . } \]
Then $W\longrightarrow U$ is faithfully flat and $u$ lifts to $v\in Z(W)$
and hence to $\phi v\in X(W)$.\\
To see the assertion about flat $f$, note first that a faithfully flat map
is trivially an epimorphism of flat sheaves. Secondly, if $f$ is flat and
an epimorphism of flat sheaves, then there is some $\phi:Z\longrightarrow X$
as in ii) and the composition $f\phi$ is surjective (on the topological
spaces underlying these affine schemes), hence so is $f$, i.e. $f$ 
is faithfully flat, \cite{Bou}, ch. II, \S2, no 5, Corollary 4,ii). The injectivity of $\psi_f$ is a special case of \cite{Bou}, I, \S3, Proposition 8, i).
\end{proof}

We have a similar result for epimorphisms of algebraic stacks.

\begin{prop}\label{epiofstack} Let $(f_0,f_1):(X_0,X_1)\longrightarrow
(Y_0,Y_1)$ be a $1$-morphism of flat Hopf algebroids inducing $f:\Xg
\longrightarrow\Yg$ on associated algebraic stacks and write
$\alpha:=t\pi_2:\cart{X_0}{f_0,Y_0,s}{Y_1}\longrightarrow Y_0$. Then the following
are equivalent:\\
i) $f$ is an epimorphism.\\
ii) $\alpha$ is an epimorphism of flat sheaves.\\
iii) There is some $\phi: Z\longrightarrow \cart{X_0}{f_0,Y_0,s}{Y_1}$ in
$\Aff$ such that $\alpha\phi$ is faithfully flat.\\
If these conditions hold then $\id_{\Modq{\Yg}}\longrightarrow f_*f^*$
is injective.
\end{prop}

\begin{proof}
The equivalence of i) and ii) is ``mise pour memoire'', the one of ii)
and iii) has been proved in Proposition \ref{epiofrep}. Assume that these
conditions hold and let $g:\Xg^{'}\longrightarrow\Xg$ be any morphism of algebraic stacks. Assume that $\id_{\Modq{\Yg}}\longrightarrow(fg)_*(fg)^*$ is injective. Then we have that the composition $\id_{\Modq{\Yg}}\longrightarrow f_*f^*\longrightarrow f_*g_*g^*f^*=(fg)_*(fg)^*$ is injective and hence so is $\id_{\Modq{\Yg}}\longrightarrow f_*f^*$. Taking $g:=P:X_0\longrightarrow\Xg$ the canonical presentation we see that we can assume that $\Xg=X_0$, in particular $f:X_0\longrightarrow\Yg$ is representable and affine (and an epimorphism). Now let $Q:Y_0\longrightarrow\Yg$ be the canonical 
presentation and form the cartesian diagram
\begin{equation}\label{diagram}
\xymatrix{ Z_0\ar[r]^{g_0}\ar[d]_-P & Y_0\ar[d]^-Q \\
 X_0 \ar[r]^-f & \Yg.}
\end{equation}

As $Q$ is $fpqc$ we have that $\id_{\Modq{\Yg}}\longrightarrow f_*f^*$ is injective if and 
only if $Q^*\longrightarrow Q^*f_*f^*\simeq g_{0,*}P^*f^*\simeq g_{0,*}g_0^*Q^*$ is injective, we used flat base change, \cite{CA} Proposition 13.1.9 and this will follow 
from the injectivity of $\id_{\Modq{Y_0}}\longrightarrow g_{0,*}g_0^*$ because $Q$ is flat.\\
As $f$ is representable and affine, $Z_0$ is an affine scheme hence, by Proposition \ref{epiofrep}, we are done because $g_0$ is an epimorphism of flat sheaves \cite{CA}, Proposition 3.8.1.
\end{proof}

There is an analogous result for monomorphisms of algebraic stacks.

\begin{prop}\label{monothm} Let $(f_0,f_1):(X_0,X_1)\longrightarrow (Y_0,Y_1)$
be a $1$-morphism of flat Hopf algebroids, $P:X_0\longrightarrow\Xg$ the rigidified algebraic stack associated with $(X_0,X_1)$, $f:\Xg\longrightarrow\Yg$ the associated
$1$-morphism of algebraic stacks, $\Theta: f^*f_*\longrightarrow\id_{\Modq{\Xg}}$
the adjunction and $\beta=(s,f_1,t):X_1\longrightarrow\cartt{X_0}{f_0,Y_0,s}{Y_1}{t,Y_0,f_0}{X_0}$. Then the following are equivalent:\\
i) $f$ is a monomorphism.\\
ii) $\beta$ is an isomorphism.\\
iii) $\Theta_{P_*{\cal O}_{X_0}}$ is an isomorphism.\\
If $f$ is representable then these conditions are equivalent to:\\
iiia) $\Theta$  is an isomorphism.\\
iiib) $f_*$ is fully faithful.
\end{prop}

\begin{remark} This result may be compared to the first assertion of Theorem 2.5
of \cite{H1}. There it is proved that $\Theta$ is an isomorphism if $f$
is a {\em flat} monomorphism.\\
In the situation of Proposition \ref{monothm}, iiib) it is natural to ask for
the essential image of $f_*$, see Proposition \ref{essimage}.\\
I do not know whether every monomorphism of algebraic stacks is representable, c.f. \cite{CA}, Corollaire 8.1.3.
\end{remark}

\begin{proof} We already know that $i)$ and $ii)$ are equivalent. Consider the 
diagram 
\[ \xymatrix{ X_0 \ar@/_/[d]_-{\Delta^{'}}\ar[r]^-P & \Xg \ar@/_/[d]_-{\Delta_f} \ar[r]^-f & \Yg\\
\pi: \Zg \ar[u]_-{\pi_1^{'}} \ar[r]_-{P^{'}} & \Xg {\times\atop f,\Yg,f} \Xg 
\ar[u]_-{\pi_1} \ar[r]_-{\pi_2} & \Xg \ar[u]_-f } \]
in which the squares made of straight arrows are cartesian. As $fP$ is representable and affine, we have $fP=\underline{\Spec}(f_*P_*\Oh_{X_0})$, c.f. \cite{CA} 14.2,
and $\pi=\underline{\Spec}(f^*f_*P_*\Oh_{X_0})$. We know that $i)$
is equivalent to the diagonal of $f$, $\Delta_f$, being an isomorphism \cite{CA}, Remarque 2.3.1. As $\Delta_f$ is a section of $\pi_1$ this is equivalent to $\pi_1$ 
being an isomorphism. As $P$ is an epimorphism, this is equivalent to 
$\pi_1^{'}$ being an isomorphism by \cite{CA}, Proposition 3.8.1. Of course, $\pi_1^{'}$
admits $\Delta^{'}:=(\id_{X_0},\Delta_f P)$ as a section so, finally, $i)$
is equivalent to $\Delta^{'}$ being an isomorphism. One checks that $\Delta^{'}=\underline{\Spec}(\Theta_{P_*\Oh_{X_0}})$ and this proves the equivalence
of $i)$ and $iii)$.\\
Now assume that $f$ is representable and a monomorphism. We will show 
that $iiia)$ holds. Consider the cartesian diagram
\[ \xymatrix{ Z \ar[r]^-{f^{'}}\ar[d]_-P & Y_0 \ar[d]^-Q\\
\Xg\ar[r]^-f & \Yg.} \]
We have
\[ P^*f^*f_* \simeq f^{'*}Q^*f_* \simeq f^{'*}f^{'}_* P^* . \]
As $P^*$ reflects isomorphism, $iiia)$ will hold if the adjunction $f^{'*}f_*^{'}\longrightarrow\id_{\Modq{Z}}$ is an isomorphism. 
As $f$ is representable, this can be checked at the stalks of $z\in Z$, and
we can replace $f^{'}$ by the induced morphism $\Spec(\Oh_{Z,z})\longrightarrow
\Spec(\Oh_{Y_0,y})$ ($y:=f^{'}(z)$) which is a monomorphism. In particular,
we have reduced the proof of $iiia)$ to the case of affine schemes, i.e. the
following assertion: If $\phi: A\longrightarrow B$ is a ring homomorphism such that
$\Spec(\phi)$ is a monomorphism, i.e. the ring homomorphism corresponding to 
the diagonal $B\otimes_A B\longrightarrow B,\; b_1\otimes b_2\mapsto b_1b_2$
is an isomorphism, then, for any $B$-module $M$, the canonical homomorphism of 
$B$-modules $M\otimes_A B\longrightarrow M$ is an isomorphism. This is however
easy:
\[ M\otimes_A B\simeq (M\otimes_B B)\otimes_A B\simeq M\otimes_B (B \otimes_A B) \simeq M\otimes_B B\simeq M,\]
and we leave it to the reader to check that the composition of these isomorphisms
is the natural map $M\otimes_A B\longrightarrow M$.\\
Finally, the proof that $iiia)$ and $iiib)$ are equivalent is a formal manipulation
with adjunctions which we leave to the reader, and trivially $iiia)$ implies
$iii)$.
\end{proof}

\begin{prop}\label{essimage} In the situation of Proposition \ref{monothm}
assume that $f$ is representable and a monomorphism, let $Q:Y_0\longrightarrow\Yg$ be the rigidified 
algebraic stack associated with $(Y_0,Y_1)$ and form the cartesian diagram 
\begin{equation}\label{wurst}
\xymatrix{ Z_0\ar[r]^{g_0}\ar[d]_-P & Y_0\ar[d]^-Q \\
\Xg\ar[r]^-f & \Yg.}
\end{equation}
Then $Z_0$ is an algebraic space and a given $\Fh\in \Modq{\Yg}$ is in the essential image of $f_*$ if and only if $Q^*\Fh$ is in the essential image of $g_{0,*}$. Consequently, $f_*$
induces an equivalence between $\Modq{\Xg}$ and the full subcategory of 
$\Modq{\Yg}$ consisting of such $\Fh$.
\end{prop}

\begin{proof}
Firstly, $Z_0$ is an algebraic space because $f$ is representable.
We know that $f_*$ is fully faithful by Proposition \ref{monothm}, $iiib)$ and need to show that the above description of its essential image is correct. If $\Fh\simeq f_*\Gh$ then $Q^*\Fh\simeq Q^*f_*\Gh\simeq g_{0,*}P^*\Gh$ so $Q^*\Fh$ lies in the essential image of $g_{0,*}$. To see the converse, extend (\ref{wurst}) to a cartesian diagram 
\[ \xymatrix{ Z_1\ar@<-.5ex>[d]\ar@<.5ex>[d]\ar[r]^-{g_1} & Y_1\ar@<-.5ex>[d]\ar@<.5ex>[d] \\
Z_0\ar[r]^{g_0}\ar[d]_-P & Y_0\ar[d]^-Q \\
\Xg\ar[r]^-f & \Yg.} \]
Note that $\Xg\simeq [ \xymatrix{Z_1\ar@<-.5ex>[r]\ar@<.5ex>[r] & Z_0 } ]$, hence
$(Z_0,Z_1)$ is a flat groupoid (in algebraic spaces) representing $\Xg$. Now let there be given $\Fh\in\Modq{\Yg}$ and $G\in\Modq{Z_0}$
with $Q^*\Fh\simeq g_{0,*}G$. We define $\sigma$ to make the following diagram commutative:
\[ \xymatrix{ s^*Q^*\Fh\ar[r]^-{can}_-{\sim}\ar[d]_-{\sim} & t^*Q^*\Fh\ar[d]_-{\sim}\\
s^*g_{0,*}G\ar[d]_-{\sim} & t^*g_{0,*}G\ar[d]_-{\sim} \\
g_{1,*}s^*G \ar[r]^-{\sim}_-{\sigma} & g_{1,*}t^*G. }\]
As $f$ is representable and a monomorphism, so is $g_1$ and thus $g_1^*g_{1,*}\stackrel{\sim}{\longrightarrow}\id_{\Modq{Z_1}}$ and $g_{1,*}$ is fully faithful by Proposition \ref{monothm},$iiia),iiib)$. We define 
$\tau$ to make the following diagram commutative:
\[ \xymatrix{ g_1^*g_{1,*}s^*G \ar[r]^-{g_1^*(\sigma)}_-{\sim} \ar[d]_-{\sim} & g_1^*g_{1,*}t^*G\ar[d]_-{\sim} \\
s^*G\ar[r]^-{\tau} & t^*G.} \]
Then $\tau$ satisfies the cocycle condition because it does so after applying 
the faithful functor $g_{1,*}$. So $\tau$ is a descent datum on $G$, and $G$ descents to $\Gh\in\Modq{\Xg}$ with
$P^*\Gh\simeq G$ and we have $Q^*f_*\Gh\simeq g_{0,*}P^*\Gh\simeq Q^*\Fh$,
hence $f_*\Gh\simeq\Fh$, i.e. $\Fh$ lies in the essential image of $f_*$ as 
was to be shown.
\end{proof}

To conclude this subsection we give the proof of Theorem \ref{isothm}
the notations and assumptions of which we now resume.

\begin{proofof} Theorem \ref{isothm}.
If iii) holds then $f$ is an epimorphism and a monomorphism by proposition
\ref{epiofstack}, $iii)\Rightarrow i)$ and Proposition \ref{monothm}, $ii)\Rightarrow i)$ hence i) holds by \cite{CA}, Corollaire 3.7.1. The proof that i) implies ii) is left to the reader
and we assume that ii) holds. Since $(f^*,f_*)$ is an adjoint pair of 
functors, $f_*$ is a quasi-inverse for $f^*$ and $\Theta: f^*f_*\longrightarrow\id_{\Modq{\Xg}}$
is an isomorphism so $\beta$ is an isomorphism by Proposition \ref{monothm}, $iii)\Rightarrow ii)$. As $f^*$
is in particular exact and faithful, $\alpha$ is faithfully flat by Proposition
\ref{flatness1}, $ii)\Rightarrow iii)$ and iii) holds.
\end{proofof}

\vspace{1.5cm}

\section{Landweber exactness and change of rings}\label{ringchange}

In this section we will use the techniques from section \ref{morphisms}
to give a short and conceptional proof of the fact that Landweber exact 
$\BP_*$-algebras of the same height have equivalent categories of comodules. In
fact, we will show that the relevant algebraic stacks are $1$-isomorphic.\\
Let $p$ be a prime number. We will study the algebraic
stack associated with the flat Hopf algebroid $(\BP_*,\BP_*\BP)$ where $\BP$
denotes Brown-Peterson homology at $p$.\\
We will work over $S:=\Spec(\Z_{(p)})$, i.e. $\Aff$ will be the category of 
$\Z_{(p)}$-algebras with its {\em fpqc} topology. We refer the reader to \cite{R1}, 
Chapter 4 for basic facts about $\BP$, e.g. $\BP_*=\Z_{(p)}[v_1,\ldots]$ where
the $v_i$ denote either the Hazewinkel- or the Araki-generators, it does not
matter but the reader is free to make a definite choice at this point if 
she feels like doing so.\\
$(V:=\Spec(\BP_*),W:=\Spec(\BP_*\BP))$ is a flat Hopf algebroid and we denote
by $P:V\longrightarrow\Xg_{FG}$ the corresponding rigidified algebraic
stack. We refer the reader to section \ref{stackformalgroups} for
an intrinsic description of the stack $\Xg_{FG}$.\\ 
For $n\geq 1$ the ideal $I_n:=(v_0,\ldots,v_{n-1})\subseteq \BP_*$ is an invariant 
prime ideal where we agree that $v_0:=p$, $I_0:=(0)$ and $I_{\infty}:=(v_0,v_1,\ldots)$.\\
As explained in subsection \ref{modules}, corresponding to these invariant ideals
there is a sequence of closed substacks
\[ \Xg_{FG}=\Zg^0\supseteq\Zg^1\supseteq\ldots\supseteq\Zg^{\infty}.\]
We denote by $\Ug^n:=\Xg_{FG}-\Zg^n$ ($0\leq n \leq
\infty$) the open substack complementary to $\Zg^n$ and have an ascending 
chain 
\[ \emptyset=\Ug^0\subseteq\Ug^1\subseteq\ldots\subseteq\Ug^{\infty}\subseteq\Xg_{FG}.\]
For $0\leq n<\infty$, $I_n$ is finitely generated, hence the open immersion
$\Ug^n\subseteq\Xg_{FG}$ is quasi-compact and $\Ug^n$ is an algebraic stack.
However, $\Ug^{\infty}$ is not algebraic: If it was, it could be covered by an 
affine (hence quasi-compact) scheme and the open covering $\Ug^{\infty}=\cup_{n\geq 0,n\neq\infty}\Ug^n$ would allow a finite subcover, which it does not.\\

\subsection{The algebraic stacks associated with Landweber exact $\BP_*$-algebras}\label{landweber}

In this subsection we prove our main result, Theorem \ref{image}, which determines
the stack theoretic image of a morphism $X_0\longrightarrow\Xg_{FG}$ corresponding
to a Landweber exact $\BP_*$-algebra. It turns out that the same arguments apply more
generally to morphisms $X_0\longrightarrow\Zg^n$ for any $n\geq 0$ 
and we work in this generality from the very beginning.\\
Fix some $0\leq n<\infty$. The stack $\Zg^n$ is associated with the flat Hopf algebroid $(V_n,W_n)$ where $V_n:=\Spec(\BP_*/I_n)$ and $W_n:=\Spec(\BP_*\BP/I_n\BP_*\BP)$, the flatness of this Hopf algebroid is established by direct inspection, and we 
have a cartesian diagram
\begin{equation}\label{diagram3}
 \xymatrix{ W_n\ar@<-.5ex>[d]\ar@<.5ex>[d]\ar@{^{(}->}[r] & W=W_0\ar@<-.5ex>[d]\ar@<.5ex>[d]\\
V_n \ar[d]_-{Q_n}\ar@{^{(}->}[r]^-{i_n} & V=V_0 \ar[d]^-Q \\
\Zg^n \ar@{^{(}->}[r] & \Xg_{FG} }
\end{equation}
in which the horizontal arrows are closed immersions.\\
We have an ascending chain of open substacks
\[ \emptyset=\Zg^n\cap\Ug^n\subseteq\Zg^n\cap\Ug^{n+1}\subseteq\ldots\subseteq\Zg^n\cap\Ug^{\infty}\subseteq\Zg^n. \]

Let $X_0\stackrel{\phi}{\longrightarrow}V_n$ be a morphism in $\Aff$ 
corresponding to a morphism of rings $\BP_*/I_n\longrightarrow R:=\Gamma(X_0,\Oh_{X_0})$. Slightly
generalising Definition 4.1 of \cite{H1} we define the height of $\phi$ as
\[ \height(\phi):=\max \{ N\geq 0 | R/I_NR \neq 0 \} \]
which may be $\infty$ and we agree to put $\height(\phi):=-1$ in case $R=0$, i.e. 
$X_0=\emptyset$. Recall that a geometric point of $X_0$ is a morphism
$\Omega\stackrel{\alpha}{\longrightarrow}X_0$ in $\Aff$ where $\Omega=\Spec(K)$
is the spectrum of an algebraically closed field $K$. The composition
$\Omega\stackrel{\alpha}{\longrightarrow}X_0\stackrel{\phi}{\longrightarrow}V_n\stackrel{i_n}{\hookrightarrow} V$ specifies a $p$-typical formal group law over 
$K$ and $\height(i_n\phi\alpha)$ is the height of this formal group law.
The relation between $\height(\phi)$ and the height of formal group laws is the following.

\begin{prop}\label{heightfibre} In the above situation we have
\[ \height(\phi)=\max\{\height(i_n\phi\alpha)|\alpha:\Omega\longrightarrow X_0\mbox{ a geometric point} \}, \]
with the convention that $\max\;\emptyset=-1$.
\end{prop}

This Proposition means that $\height(\phi)$ is the maximum height in a geometric
fibre of the formal group law over $X_0$ parametrised by $i_n\phi$.

\begin{proof} Clearly, $\height(i_n\phi\psi)\leq\height(\phi)$ for any morphism 
$\psi:Y\longrightarrow X_0$ in $\Aff$. For any $0\leq N^{'}\leq\height(\phi)$
we have that $I_{N^{'}}R\neq R$ so there is a maximal ideal of $R$
containing $I_{N^{'}}R$ and a geometric point $\alpha$ of $X_0$ supported at
this maximal ideal will satisfy $\height(i_n\phi\alpha)\geq N^{'}$.
\end{proof}

Another geometric interpretation of $\height(\phi)$ is given by considering the composition 
$f:X_0\stackrel{\phi}{\longrightarrow}V_n\stackrel{Q_n}{\longrightarrow}\Zg^n$.

\begin{prop}\label{heightfactor} In this situation we have
\[ \height(\phi)+1=\min\{ N\geq 0|f \mbox{ factors through } \Zg^n\cap\Ug^N\hookrightarrow\Zg^n \} \]
with the convention that $\min\;\emptyset=\infty$ and $\infty+1=\infty$.
\end{prop}

\begin{proof} For any $\infty>N\geq n$ we have a cartesian square 
\begin{equation}\label{diagram2}
\xymatrix{ V_n^N\ar[r]^-j\ar[d] & V_n\ar[d]^-{Q_n} \\
\Zg^n\cap\Ug^N\ar[r]^-i & \Zg^n}
\end{equation}

where $V_n^N=V_n-\Spec(\BP_*/I_N)=\bigcup_{i=n}^{N-1}\Spec((\BP_*/I_n)[v_i^{-1}])$
hence $f$ factors through $i$ if and only if $\phi:X_0\longrightarrow V_n$
factors through $j$. As $j$ is an open immersion, this is equivalent to
$|\phi|(|X_0|)\subseteq |V_n^N| \subseteq |V_n|$ where $| \cdot |$ denotes 
the topological space underlying a scheme. But this condition can 
be checked using geometric points and the rest is easy, using Proposition \ref{heightfibre}.
\end{proof}

Recall from \cite{H1}, Definition 2.1 that, if $(A,\Gamma)$ is a flat Hopf algebroid, 
an $A$-algebra $f: A\longrightarrow B$ is said to be {\em Landweber exact} over $(A,\Gamma)$ if the functor $M\mapsto M\otimes_A B$ from $\Gamma$-comodules
to $B$-modules is exact. For $(X_0:=\Spec(A),X_1:=\Spec(\Gamma))$, 
$\phi:=\Spec(f): Y_0:=\Spec(B)\longrightarrow X_0$ and $P: X_0\longrightarrow\Xg$
the rigidified algebraic stack associated with $(X_0,X_1)$ this exactness is equivalent to the flatness of the
composition $Y_0\stackrel{\phi}{\longrightarrow}X_0\stackrel{P}{\longrightarrow}\Xg$ because the following square of functors commutes up to
natural isomorphism
\[ \xymatrix{ (P\phi)^*:\Modq{\Xg} \ar[r] \ar[d]^{\simeq} & \Modq{Y_0} \ar[d]^{\simeq} \\
\Gamma\mbox{-comodules} \ar[r]^{M\mapsto M\otimes_A B} & B\mbox{-modules,}}  \]
where the horizontal equivalences are those given by (\ref{boing}).\\
In case $\Xg=\Zg^n$
this flatness has the following decisive consequence which paraphrases
the fact that the image of a flat morphism is stable under generalisation.

\begin{prop}\label{allheights} Assume that $n\geq 0$ and that $\phi:\emptyset\neq X_0\longrightarrow V_n$
is Landweber exact of height $N:=\height(\phi)$ (hence $n\leq N\leq\infty$). Then 
for any $n\leq j\leq N$ there is a geometric point $\alpha:\Omega\longrightarrow
X_0$ such that $\height(i_n\phi\alpha)=j$.
\end{prop}

\begin{proof} Let $\phi$ correspond to $\BP_*/I_n\longrightarrow R$. We 
first note that $v_n,v_{n+1},\ldots\in R$ is a regular sequence by Proposition \ref{landweberflat} below. Now assume that $N<\infty$ and fix $n\leq j\leq N$. Then $v_j\in R/I_{j-1}R\neq 0$ is not a zero divisor
and thus there is a minimal prime ideal of $R/I_{j-1}R$ not containing $v_j$.
A geometric point supported at this prime ideal solves the problem. In the
remaining case $j=N=\infty$ we have $R/I_{\infty}R\neq 0$ and any geometric point of this 
ring solves the problem.
\end{proof}

The main result of this paper is the following.

\begin{theorem}\label{image}
Assume that $n\geq 0$ and that $\emptyset\neq X_0\longrightarrow V_n$ is Landweber exact of height $N$ (hence 
$n\leq N\leq\infty)$. Let $(X_0,X_1)$ be the Hopf algebroid induced 
from $(V,W)$ by the composition $X_0\stackrel{\phi}{\longrightarrow}V_n
\stackrel{i_n}{\hookrightarrow} V$. Then $(X_0,X_1)$ is a flat Hopf algebroid and its 
associated algebraic stack is given as
\[ [ \xymatrix{ X_1\ar@<-.5ex>[r]\ar@<.5ex>[r] & X_0 } ] \simeq \Zg^n\cap\Ug^{N+1}
\mbox{ if }N\neq\infty\mbox{ and} \]
\[ [ \xymatrix{ X_1\ar@<-.5ex>[r]\ar@<.5ex>[r] & X_0 } ] \simeq \Zg^n\mbox{ if }N=\infty. \]
\end{theorem}

\begin{proof}
Note that $(X_0,X_1)$ is also induced from the flat Hopf algebroid $(V_n,W_n)$
along $\phi$ and thus is a flat Hopf algebroid using the final statement of Proposition \ref{flatness2} and the 
Landweber exactness of $\phi$. We first assume that $N\neq\infty$. 
Then by Proposition \ref{heightfactor} the composition $X_0\stackrel{\phi}{\longrightarrow}V_n\longrightarrow\Zg^n$ 
factors as $X_0\stackrel{\psi}{\longrightarrow}\Zg^n\cap\Ug^{N+1}\stackrel{i}{\longrightarrow}\Zg^n$ and $\psi$ is flat because $i$ is an open immersion and 
$X_0\longrightarrow\Zg^n$ is flat by assumption. By Proposition \ref{inducediso}
we will be done if we can show that $\psi$ is in fact faithfully flat. For this
we consider the presentation $\Zg^n\cap\Ug^{N+1}\simeq [ \xymatrix{ W_n^{N+1}\ar@<-.5ex>[r]\ar@<.5ex>[r] & V_n^{N+1} } ]$ given by the cartesian diagram
\[ \xymatrix{ W_n^{N+1} \ar[r]\ar@<-.5ex>[d]\ar@<.5ex>[d] & W_n \ar@<-.5ex>[d]\ar@<.5ex>[d] \\
V_n^{N+1} \ar[r]\ar[d] & V_n \ar[d]^-{Q_n} \\
\Zg^n\cap\Ug^{N+1} \ar[r] & \Zg^n} \]
and note that $\psi$ lifts to $\rho: X_0\longrightarrow V_n^{N+1}$ and induces
$\alpha:=t\pi_2: \cart{X_0}{\rho,V_n^{N+1},s}{W_n^{N+1}}\longrightarrow V_n^{N+1}$ which
is flat and we need it to be faithfully flat to apply Proposition \ref{flatness1}, $iii)\Rightarrow iv)$ and conclude that $\psi$ is faithfully flat. So we have to prove that $\alpha$ is surjective on the topological
spaces underlying the schemes involved.\\
This surjectivity can be checked on geometric points and for
any such geometric point $\Omega\stackrel{\mu}{\longrightarrow}V_n^{N+1}$ we have that 
$j:=\height(\Omega\stackrel{\mu}{\longrightarrow}V_n^{N+1}\longrightarrow V_n\stackrel{i_n}{\hookrightarrow}V)$
satisfies $n\leq j\leq N$. By Proposition \ref{allheights} there is a geometric
point $\Omega^{'}\stackrel{\nu}{\longrightarrow}X_0$ with $\height(\Omega^{'}\stackrel{\nu}{\longrightarrow}X_0\longrightarrow V_n\stackrel{i_n}{\hookrightarrow}V)=j$ and we can assume 
that $\Omega=\Omega^{'}$ because the corresponding fields have the same
characteristic namely $0$ if $j=0$ and $p$ otherwise. As any two
formal group laws over an algebraically closed field having the same
height are isomorphic we find some $\sigma:\Omega\longrightarrow W_n^{N+1}$
fitting into a commutative diagram
\[ \xymatrix{ X_0{\times\atop \rho,V_n^{N+1},s} W_n^{N+1}\ar[r]^-{\alpha} & V_n^{N+1}\\
\Omega.\ar[u]^-{(\nu,\sigma)}\ar[ur]_-{\mu} & } \]
As $\mu$ was arbitrary this shows that $\alpha$ is surjective.
We leave the obvious modifications for the case $N=\infty$ to the reader.
\end{proof}

To conclude this subsection we explain the relation of Landweber exactness and Landweber's
regularity condition. This is well-known to the expert and in fact has been worked out in detail
in \cite{franke}, section 3, Theorem 8 but we include it here anyway. Fix some $n\ge 0$ and let 
$\phi:\BP_*/I_n \longrightarrow R$ be a $\BP_* / I_n$-algebra. Then Landweber's condition is
\begin{equation} \label{landreg}
\mbox{ The sequence }\phi(v_n),\phi(v_{n+1}),\ldots\,\in R\mbox{ is regular.}
\end{equation}

\begin{prop}\label{landweberflat}
In the above situation, (\ref{landreg}) holds if and only if the composition $\Spec(R)\longrightarrow\Spec(\BP/I_n)\longrightarrow\Zg^n$ is flat.
\end{prop}

\begin{proof} From \cite{millerrav}, Proposition 2.2 we know that the restriction of 
$f^*: \Modq{\Zg^n}\longrightarrow\Modq{\Spec(R)}$ to finitely presented comodules is exact
if and only if (\ref{landreg}) holds. But $f^*$ itself is exact, and hence $f$ is flat, if 
and only if its above restriction is exact because any $\BP_*\BP/I_n$-comodule is the filtering
direct limit of finitely presented comodules. This was pointed out to me by N. Strickland. In case $n=0$ this is \cite{millerrav}, Lemma 2.11
and the general case follows from \cite{H2}, Proposition 1.4.1,e), Proposition 1.4.4, Lemma
1.4.6 and Proposition 1.4.8.
\end{proof}

\subsection{Equivalence of comodule categories and change of rings}

In this subsection we will spell out some consequences of the above results in the language of comodules but need some elementary preliminaries first.\\
Let $A$ be a ring, $I=(f_1,\ldots, f_n)\subseteq A$ ($n\geq 1$) a finitely
generated ideal and $M$ an $A$-module. We have a canonical map
\[ \bigoplus_i M_{f_i}\longrightarrow\bigoplus_{i<j}M_{f_if_j},\; (x_i)_i\mapsto \left( \frac{x_i}{1}-\frac{x_j}{1} \right)_{i,j} \]
and a canonical map \[ \alpha_M: M\longrightarrow\ker(\bigoplus_i M_{f_i}\longrightarrow\bigoplus_{i<j}M_{f_if_j}).\]
For $X:=\Spec(A)$, $Z:=\Spec(A/I)$, $j: U:=X-Z\hookrightarrow X$ the open immersion 
and $\Fh$ the quasi-coherent $\Oh_X$-module corresponding to $M$, $\alpha_M$
corresponds to the adjunction $\Fh\longrightarrow j_*j^*\Fh$.
Note that $\ker(\alpha_M)$ is the $I$-torsion submodule of $M$.
The cokernel of $\alpha_M$ corresponds to the local cohomology $H^1_Z(X,\Fh)$, 
c.f. \cite{local}. We say that $M$ is $I$-local if $\alpha_M$ is an
isomorphism. A quasi-coherent $\Oh_X$-module $\Fh$ is in the essential image 
of $j_*$ if and only if $\Fh\longrightarrow j_*j^*\Fh$ is an isomorphism
if and only if the $A$-module corresponding to $\Fh$ is $I$-local. If $n=1$
then $M$ is $I=(f_1)$-local if and only if $f_1$ acts invertibly on $M$.\\
We now formulate a special case of Proposition \ref{essimage} in terms
of comodules.

\begin{prop}\label{subcomod} For any $n\geq 0$ the category $\Modq{\Zg^n}$
is equivalent to the full subcategory of $\BP_*\BP$-comodules $M$ such that $I_nM=0$.\\
For any $0\leq n\leq N<\infty$ the category $\Modq{\Zg^n\cap\Ug^{N+1}}$ is equivalent to the full subcategory of $\BP_*\BP$-comodules $M$ such that $I_nM=0$ and
$M$ is $I_{N+1}/I_n$-local as a $\BP_*/I_n$-module.
\end{prop}

\begin{remark} We know from (\ref{boing}) that $\Modq{\Zg^n}$ is equivalent to the category of
$\BP_*\BP/I_n$-comodules. The alert reader will have noticed that we have not
yet mentioned any graded comodules. This is not sloppy terminology, we really mean 
comodules without any grading even though the flat Hopf algebroids are all graded. 
However, it is easy to take the grading into account, in particular all results 
of this subsection have analogues for graded comodules, c.f. Remark \ref{grading}.
\end{remark}

\begin{proof} Fix $0\leq n<\infty$. The $1$-morphism $\Zg^n\hookrightarrow\Xg_{FG}$ is representable and a closed immersion (in particular a monomorphism) because its base change along $V\longrightarrow\Xg_{FG}$ is a closed immersion and being a closed immersion is $fpqc$-local on the base, \cite{ega42}, 2.7.1, $xii)$. Proposition \ref{essimage}
identifies $\Modq{\Zg^n}$ with the full subcategory of $\Modq{\Xg_{FG}}$ 
consisting of those $\Fh$ such that $Q^*\Fh\simeq i_{n,*}G$ for some 
$G\in\Modq{V_n}$ ( with notations as in (\ref{diagram3})). Identifying, as in subsection \ref{modules}, $\Modq{\Xg_{FG}}$
with the category of $\BP_*\BP$-comodules, $\Fh$ corresponds to some $\BP_*\BP$-comodule
$M$ and $Q^*\Fh$ corresponds to the $\BP_*$-module underlying $M$.
So the condition of Proposition \ref{essimage} is that the $\BP_*$-module
$M$ is in the essential image of $i_{n,*}$, i.e. $M$ is an $\BP_*/I_n$-module,
i.e. $I_nM=0$.\\
Now fix $0\leq n \leq N<\infty$. We apply Proposition \ref{essimage} 
to $i:\Zg^n\cap\Ug^{N+1}\longrightarrow\Xg_{FG}$  which is representable and a quasi-compact immersion (in particular a monomorphism) because it 
sits in a cartesian diagram
\[ \xymatrix{ V_n^{N+1} \ar[d]\ar[r]^-j & V \ar[d]^-Q \\
\Zg^n\cap\Ug^{N+1} \ar[r]^-i & \Xg_{FG}, } \]
c.f. (\ref{diagram2}), in which $j$ is a quasi-compact immersion and one uses
\cite{ega42}, 2.7.1, $xi)$ as above. 
Arguing as above, we are left with identifying the essential image 
of $j_*$ which, as explained at the beginning of this subsection, corresponds to the $\BP_*$-modules $M$
such that $I_nM=0$ and $M$ is $I_{N+1}/I_n$-local as a $\BP_*/I_n$-module.
\end{proof}

\begin{cor}\label{comodcat}
Let $n\geq 0$ and let $\BP_*/I_n\longrightarrow R\neq 0$ be Landweber exact of height $N$ (hence $n\leq N \leq\infty$). Then $(R,\Gamma):=(R,R\otimes_{\BP_*}\BP_*\BP\otimes_{\BP_*}R)$ is a flat Hopf algebroid and its
category of comodules is equivalent to the full subcategory of $\BP_*\BP$-comodules
$M$ such that $I_nM=0$ and $M$ is $I_{N+1}/I_n$-local as a $\BP_*/I_n$-module. The last condition is to be ignored in case $N=\infty$. 
\end{cor}
\begin{proof}
By Theorem \ref{image}, $(R,\Gamma)$ is a flat Hopf algebroid with associated
algebraic stack $\Zg^n\cap\Ug^{N+1}$ (resp. $\Zg^n$ if $N=\infty$). So the category
of $(R,\Gamma)$-comodules is equivalent to $\Modq{\Zg^n\cap\Ug^{N+1}}$ (resp. $\Modq{\Zg^n}$). 
Now use Proposition \ref{subcomod}.
\end{proof}

\begin{remark} The case $n=0$ of Corollary \ref{comodcat} corresponds to the situation treated in \cite{H1} where, translated into the present terminology, $\Modq{\Ug^{N+1}}$ is identified as a {\em localisation} of $\Modq{\Xg_{FG}}$. 
This can be done because $f:\Ug^{N+1}\longrightarrow\Xg_{FG}$ is flat, hence $f^*$ exact.
To relate more generally $\Modq{\Zg^n\cap\Ug^{N+1}}$ to $\Modq{\Xg_{FG}}$
it seems more appropriate to identify the former as a full subcategory
of the latter as we did above. However, using Proposition 1.4 of {\em loc. cit.}
and Proposition \ref{monothm} one sees that $\Modq{\Zg^n\cap\Ug^{N+1}}$
is equivalent to the localisation of $\Modq{\Xg_{FG}}$ with respect to 
all morphisms $\alpha$ such that $f^*(\alpha)$ is an isomorphism where
$f:\Zg^n\cap\Ug^{N+1}\longrightarrow\Xg_{FG}$ is the immersion. As $f$ is not flat
for $n\geq 1$ this condition seems less tractable than the one in Corollary \ref{comodcat}.
\end{remark}

Of course, equivalences of comodule categories give rise to change of rings
theorems and we refer to \cite{H1} for numerous examples (in the case $n=0$)
and only point out the following, c.f. \cite{R2}, Theorem B.8.8 for the 
notation and a special case: If $n\geq 1$ and $M$ is a $\BP_*\BP$-comodule such that 
$I_nM=0$ and $v_n$ acts invertibly on $M$ then
\[ \Ext^*_{\BP_*\BP}(\BP_*,M)\simeq\Ext^*_{\Sigma(n)}(\F_p[v_n,v_n^{-1}],M\otimes_{\BP_*}
\F_p[v_n,v_n^{-1}]) . \]

In fact, this is clear from the case $n=N$ of Corollary \ref{comodcat}
applied to the obvious map $\BP_*/I_n\longrightarrow\F_p[v_n,v_n^{-1}]$ which is Landweber exact of height $n$.\\
To make a final point, in \cite{H1} we also find many of the fundamental results of
\cite{L} generalised to Landweber exact algebras whose induced Hopf algebroids are presentations of our $\Ug^{N+1}$. One may generalise these results further to the present case, i.e. to $\Zg^n\cap\Ug^{N+1}$ for $n\geq 1$,
but again we leave this to the reader and only point out an example: In the situation of Corollary \ref{comodcat} every non-zero graded
$(R,\Gamma)$-comodule has a non-zero primitive.\\
To prove this, consider the comodule as a quasi-coherent sheaf $\Fh$ on $\Zg^n\cap\Ug^{N+1}$ and use
that the primitives we are looking at are $H^0(\Zg^n\cap\Ug^{N+1},\Fh)\simeq
H^0(\Xg_{FG},f_*\Fh)\neq 0$ because $f_*$ is faithful and using the 
result of P. Landweber that every non-zero graded $\BP_*\BP$-comodule
has a non-zero primitive.

\section{The stack of formal groups}\label{stackformalgroups}

In this section we take a closer look at the algebraic stacks associated with the
flat Hopf algebroids $(\MU_*,\MU_*\MU)$ and $(\BP_*,\BP_*\BP)$.\\
A priori, these stacks are given by the abstract procedure of stackification and in many instances one can work with this definition directly, the results of the previous sections are an example of this. 
For future investigations, e.g. those initiated in \cite{G}, it might be
useful to have the genuinely geometric description of these stacks which we
propose to establish in this section.\\
For this, we require a good notion of formal scheme over an arbitrary affine
base as given by N. Strickland \cite{formalgrp} and we quickly recall some of his results now.\\
The category $X_{fs,\Z}$ of formal schemes over $\Spec(\Z)$ is defined to be
the ind-category of $\Aff_{\Z}$ which we consider as usual as a full
subcategory of the functor category $C:=\underline{\Hom}(\Aff_{\Z}^{op},\mathrm{Sets})$, c.f. \cite{formalgrp}, Definition 4.1 and \cite{SGA4}, expos\'e I, 8. A formal ring is by definition
a linearly topologised Hausdorff and complete ring and $\FRings$ denotes the category
of formal rings with continuous ring homomorphisms. Any ring can be considered as 
a formal ring by giving it the discrete topology. There is a fully faithful functor $\Spf: \FRings^{\op}\longrightarrow X_{fg,\Z}\subset C$ \cite{formalgrp}, section 4.2
given by 
\[ \Spf(R)(S):=\Hom_{\FRings}(R,S)=\mbox{colim}_I \Hom_{\Rings}(R/I,S),\]
the limit being taken over the directed set of open ideals $I\subseteq R$.\\
In particular, any ring $R$ can be considered as a formal scheme over $\Z$ 
and we thus get the category $X_{fs,R}:=X_{fs,\Z}/\Spf(R)$ of formal schemes 
over $R$. For varying $R$, these categories assemble into an $fpqc$-stack $X_{fs}$ over $\Spec(\Z)$ which we call the stack of formal schemes \cite{formalgrp}, Remark 2.58, Proposition 4.51 and  Remark 4.52.\\
Define $X_{fgr}$ to be the category of commutative group objects in $X_{fs}$.
Then $X_{fgr}$ is canonically fibred over $\Aff_{\Z}$ and is in fact an
$fpqc$-stack over $\Spec(\Z)$ because being a commutative group object can be expressed by the
existence of suitable structure morphisms making appropriate diagrams commute.
Finally, define $X\subseteq X_{fgr}$ to be the substack of those objects which are $fpqc$-locally isomorphic to $(\Aone,0)$ as
{\em pointed formal schemes} (of course, a formal group is considered as a pointed formal schemes via its zero section). It is clear that $X\subseteq X_{fgr}$
is in fact a substack and in particular is itself an $fpqc$-stack over
$\Spec(\Z)$ which we will call the stack of formal groups. We will see
in a minute that $X$ (unlike $X_{fgr}$) is in fact an algebraic stack.\\
Our first task will be to determine what formal schemes occur in the fibre category $X_R$ for a given ring $R$. This requires some notation:\\
For a locally free $R$-module $V$ of rank one we denote by $\hat{SV}$ the symmetric algebra of $V$ over $R$ completed with respect to its augmentation ideal.
This $\hat{SV}$ is a formal ring. The diagonal morphism $V\longrightarrow V\oplus V$ induces a structure of
formal group on $\Spf(\hat{SV})$. Indeed, for any faithfully flat extension
$R\longrightarrow R^{'}$ with $V\otimes_R R^{'}\simeq R'$ we have
$\Spf(\hat{SV})\times_{\Spec(R)} \Spec(R^{'})\simeq\hat{\Ge}_{a,R^{'}}$
in $X_{R^{'}}$. On the other hand, denote by $\Sigma(R)$ the set of 
isomorphism classes of pointed formal schemes in $X_R$. We have a map
$\rho_R:\Pic(R)\longrightarrow\Sigma(R)\;,\; [V]\mapsto [\Spf(\hat{SV})]$.

\begin{prop}\label{forms}

For any ring $R$, the map $\rho_R:\Pic(R)\longrightarrow\Sigma(R)$ is bijective.
\end{prop}

\begin{remark} For suitable rings $R$ we can compare the above construction of the category of formal groups over $R$ to more traditional ones: 
If $R$ is pseudocompact and local \cite{conrad}, Definition 1.1.4 then, using Proposition
\ref{forms}, one can check that $X_R$ is the groupoid of one dimensional commutative formal Lie groups over $R$ in the sense of \cite{conrad}, Definition 3.3.2.
\end{remark}

\begin{proofof} Proposition \ref{forms}. By definition, $\Sigma(R)$ is the set 
of $fpqc$-forms of the pointed formal scheme $(\Aone,0)$ over $R$. We thus have a
$\check{\mathrm{C}}$ech-cohomological description
\[ \Sigma(R)\simeq\CH^1(R,\underline{\Aut}(\Aone,0))=\mbox{colim}_{R\longrightarrow R^{'}} \CH^1(R^{'}/R,\underline{\Aut}(\Aone,0)),\]
where $G^0:=\underline{\Aut}(\Aone,0)$ is the sheaf of automorphisms of the pointed 
formal scheme $(\Aone,0)$ over $R$ and the limit is taken over all faithfully flat extensions $R\longrightarrow R^{'}$. For an arbitrary $R$-algebra $R^{'}$ we can
identify \[ G^{0}(R^{'})= \{ f\in R^{'}[[t]] \; | \; f(0)=0,\; f^{'}(0)\in R^* \} \]
with the multiplication of the right hand side being substitution of
power series. We have a split epimorphism $\pi:G^{0}\longrightarrow\Ge_m$ given
on points by $\pi(f):=f^{'}(0)$ with kernel $G^{1}:=\ker(\pi)$ and we define
more generally for any $n\geq 1$, $G^n(R^{'}):= \{ f \in G^{0}(R^{'}) \; | \; f=1+ O(t^n) \}$. For any $n\geq 1$ we have an epimorphism $G^n\longrightarrow\Ge_a,\;
f=1+\alpha t^n+O(t^{n+1})\mapsto\alpha$ with kernel $G^{n+1}$.
One checks that the $G^n$ are a descending chain of normal subgroups in $G^0$ defining for every $R$-algebra $R^{'}$ a structure of complete Hausdorff topological group on $G^0(R^{'})$.\\
Using $\CH^1(R^{'}/R,\Ge_a)=0$ and an approximation argument shows that
$\CH^1(R^{'}/R,G^1)=0$ for any $R$-algebra $R^{'}$, hence the map $\phi:\CH^1(R,G^0)\longrightarrow\CH^1(R,\Ge_m)$ induced by $\pi$ is injective, and as $\pi$ is split we see that 
$\phi$ is a bijection. As $\CH^1(R,\Ge_m)\simeq\Pic(R)$ we have obtained a bijection $\Sigma(R)\simeq\Pic(R)$ and unwinding the definitions shows that it coincides with $\rho_R$.
\end{proofof}

The stack $X$ carries a canonical line bundle:\\
For any ring $R$ and $G\in X_R$ we can construct the locally free rank one $R$-module $\omega_{G/R}$ as usual \cite{formalgrp}, Definition 7.1 and as its formation is
compatible with base change it defines a line bundle $\omega$ on $X$. We remark
without proof that $\Pic(X)\simeq\Z$, generated by the class of $\omega$.\\
We define a $\Ge_m$-torsor $\pi: \Xg:=\underline{\Spec}(\oplus_{\nu\in\Z}\omega^{\otimes \nu})\longrightarrow X$ \cite{CA}, 14.2 and now check that $\Xg$
is the algebraic stack associated with the flat Hopf algebroid $(\MU_*,\MU_*\MU)$.\\
For any ring $R$, the category $\Xg_R$ is the groupoid of pairs $(G/R,\; \omega_{G/R}\stackrel{\simeq}{\longrightarrow} R)$ consisting of a formal group $G/R$
together with a trivialization of the $R$-module $\omega_{G/R}$. The morphisms
in $\Xg_R$ are the isomorphisms of formal groups which respect the trivializations
in an obvious sense. Since $\omega_{\Spf(\hat{SV})/R}\simeq V$ we see from Proposition \ref{forms} that 
any $G\in\Xg_R$ is isomorphic to $(\Aone,0)$ as a pointed formal scheme over $R$.
This easily implies that the diagonal of $\Xg$ is representable and affine.
Now recall the affine scheme FGL$\simeq\Spec(\MU_*)$ \cite{formalgrp}, Example 2.6
parametrising formal group laws. We define $f:$FGL$\longrightarrow\Xg$ by 
specifying the corresponding object of $\Xg_{\FGL}$ as follows: 
We take $G:=\Aone_{\FGL}=\Spf(\MU_*[[x]])$ with the group structure induced by
a fixed choice of universal formal group law over $\MU_*$ together with
the trivialization $\omega_{G/\MU_*}=(x)/(x^2)\stackrel{\simeq}{\longrightarrow}
\MU_*$ determined by $x\mapsto 1$. We then claim that $f$ is faithfully flat and
thus $\Xg$ is an algebraic stack with presentation $f$ (this will also
imply that $X$ is an algebraic stack):\\
Given any $1$-morphism $\Spec(R)\longrightarrow\Xg$ we can assume that that 
the corresponding object of $\Xg_R$ is given as $(\Aone_R,\; (x)/(x^2)\stackrel{\simeq}{\longrightarrow} R,\; x\mapsto u)$ with the group structure on $(\Aone_R,0)$ defined by some formal group law over
$R$ and with some unit $u\in R^*$. Then $\Spec(R)\times_{\Xg} \FGL$ parametrises isomorphisms of 
formal group laws with leading term $u$. This is well-known to be representable
by a polynomial ring over $R$, hence it is faithfully flat.\\
The same argument shows that $\FGL\times_{\Xg}\FGL\simeq\FGL\times_{\Spec(\Z)}\; \SI\simeq\Spec(\MU_*\MU)$ where $\SI$ parametrises strict isomorphisms of formal group laws
\cite{R1}, Appendix A 2.1.4 and this establishes the first half of the following
result.

\begin{theorem}\label{mubp} 1) $\Xg$ is the algebraic stack associated with the flat Hopf algebroid $(\MU_*,\MU_*\MU)$.\\
2) For any prime $p$, $\Xg\times_{\Spec(\Z)}\Spec(\Z_{(p)})$ is the algebraic stack associated with the flat Hopf algebroid $(\BP_*,\BP_*\BP)$.
\end{theorem}

\begin{proof} The proof of 2) is identical to the proof of 1) given above except
that to see that the obvious $1$-morphism $\Spec(\BP_*)\longrightarrow\Xg\times_{\Spec(\Z)}\Spec(\Z_{(p)})$ is faithfully flat one has to use Cartier's theorem
saying that any formal group law over a $\Z_{(p)}$-algebra is strictly isomorphic
to a $p$-typical one , see for example \cite{R1}, Appendix A 2.1.18. 
\end{proof}

\begin{remark}\label{grading} 1) We explain how the grading of $\MU_*$ fits into the above result. The stack $\Xg$ carries a $\Ge_m$-action given on points by
\[ \alpha\cdot(G/R,\; \phi:\omega_{G/R}\stackrel{\simeq}{\longrightarrow} R):=
(G/R,\; \phi:\omega_{G/R}\stackrel{\simeq}{\longrightarrow} R\stackrel{\cdot\alpha}{\longrightarrow} R)\mbox{ for }\alpha\in R^* .\]
This action can be lifted to the Hopf algebroid $(\FGL,\FGL\times\SI)$ as in
\cite{formalgrp}, Example 2.97 and thus determines a grading of the flat Hopf algebroid
$(\MU_*,\MU_*\MU)$. As observed in {\em loc.cit} this is the usual (topological) grading except that all degrees are divided by 2.\\
2) For any $n\geq 0$ we know from section \ref{stacksandhopf} and Theorem \ref{mubp}, 1) that 
\[ \Ext^n_{\MU_*\MU}(\MU_*,\MU_*) \simeq \HH^{n}(\Xg,\Oh_{\Xg}).\]
As $\pi:\Xg\longrightarrow X$ is affine its Leray spectral sequence collapses to
an isomorphism $\HH^n(\Xg,\Oh_{\Xg})\simeq\HH^n(X,\pi_* \Oh_{\Xg})\simeq\oplus_{k \in\Z}\HH^n(X,\omega^{\otimes k})$. The comparison of gradings given in 1) implies that this isomorphism restricts, for every $k\in\Z$, to an isomorphism
\[ \Ext^{n,2k}_{\MU_*\MU}(\MU_*,\MU_*)\simeq\HH^n(X,\omega^{\otimes k}).\]
In particular, we have $\HH^*(X,\omega^{\otimes k})=0$ for all $k<0$.\\
3) As $\pi:\Xg\longrightarrow X$ is fpqc, the pull back $\pi^*$
establishes an equivalence between $\Modq{X}$ and the category of quasi-coherent $\Oh_{\Xg}$-modules equipped with a descent datum with respect to $\pi$, c.f.
the beginning of subsection \ref{modules}. One checks that a descent datum on a given $\Fh\in\Modq{\Xg}$ with respect to $\pi$ is the same as a $\Ge_m$-action 
on $\Fh$ compatible with the action on $\Xg$ given in 1). Hence $\pi^*$ gives an
equivalence between $\Modq{X}$ and the category of {\em evenly graded} $\MU_*\MU$-comodules.\\
4) The referee suggest a different way of looking at 3): Since $\Xg\longrightarrow X$ is a 
$\Ge_m$-torsor it is in particular $fpqc$ and hence the composition $\Spec(\MU_*)\longrightarrow
\Xg\longrightarrow X$ is a presentation of $X$ and one checks that the corresponding flat 
Hopf algebroid is $(\MU_*,\MU_*\MU[u^{\pm 1}])$ thereby justifying our ad hoc definition of 
$X$ in section \ref{prelim}. This again shows that $\Modq{X}$ is equivalent to the category
of evenly graded $\MU_*\MU$-comodules, this time the grading being accounted for by the
coaction of $u$.\\
5) The analogues of 1)-4) above with $X$ (resp. $\MU$) replaced by
$X\times_{\Spec(\Z)}\Spec(\Z_{(p)})$ (resp. $\BP$) hold true. 
\end{remark}

The last issue we would like to address is the stratification of $X$ by the
height of formal groups.\\
For every prime $p$ we put $Z_p^1:=X\times_{\Spec(\Z)} \Spec(\F_p)\subseteq X$.\\
The universal formal group $G$ over $Z_p^1$ comes equipped with a relative 
Frobenius $F:G\longrightarrow G^{(p)}$ which can be iterated to $F^{(h)}:G\longrightarrow G^{(p^h)}$ for all $h\geq 1$.\\
For $h\geq 1$ we define $Z^h_p \subseteq Z^1_p$ to be the locus over which the $p$-multiplication of $G$ factors through $F^{(h)}$. Clearly, $Z_p^h\subseteq X$ is a closed substack, hence 
$Z^h_p$ is the stack of formal groups over $\Spec(\F_p)$ which have height at least $h$. The stacks labeled $\Zg^n$ ($n\geq 1$) in section \ref{ringchange}
are the preimages of $Z_p^n$ under $\pi\times\id:\Xg\times_{\Spec(\Z)}\Spec(\Z_{(p)})\longrightarrow X\times_{\Spec(\Z)}\Spec(\Z_{(p)})$.\\
For any $n\geq 1$ we define the (non-closed) substack $Z^n:=\bigcup_{p\mbox{\small prime}}Z_p^n\subseteq X$ with complement $U^n:=X-Z^n$.\\
If $\MU_*\longrightarrow B$ is a Landweber exact $\MU_*$-algebra which has height $n\geq 1$ at every prime as in \cite{H1}, section 7 then the stack theoretic image of $\Spec(B)\longrightarrow\Spec(\MU_*)\longrightarrow\Xg$ is the preimage of $U^n$ under $\pi: \Xg\longrightarrow X$ which we will
write as $\Ug^n:=\pi^{-1}(U^n)\subseteq\Xg$. This can be checked as in section \ref{ringchange}
and shows that the equivalences of comodule categories proved in {\em loc. cit.} are again a consequence of the fact that the relevant algebraic stacks are 
$1$-isomorphic. We leave the details to the reader.\\
To conclude we would like to point out the following curiosity:\\
As complex $\K$-theory is Landweber exact of height $1$ over $\MU_*$ we know 
that the flat Hopf algebroid $(\K_*,\K_*\K)$ has $\Ug^1$ as its associated 
algebraic stack.
So J. Adams' computation of $\Ext^1_{\K_*\K}(\K_*,\K_*)$ implies that for any integer
$k\geq 2$ we have 
\[ |\HH^1(U^1,\omega^{\otimes k})|=2\cdot\mbox{denominator}\;(\zeta(1-k)),\]
where $\zeta$ is the Riemann zeta function and we declare the denominator of $0$ to be $1$. To check this one uses Remark \ref{grading}, 2) with $X$ replaced by $U^1$, \cite{switzer}, Proposition 19.22  and \cite{neukirch}, VII, Theorem 1.8.\\
Unfortunately, the orders of
the (known) groups $\HH^2(U^1,\omega^{\otimes k})$ have nothing to do 
with the nominators of Bernoulli-numbers.\\

\end{document}